\documentclass[a4paper,12pt]{article}
\usepackage[utf8]{inputenc}
\usepackage[backend=biber, style=alphabetic,giveninits=true, url=false, eprint=false, isbn=false ,maxbibnames=9]{biblatex}

\usepackage{amsmath}
\usepackage{amssymb}
\usepackage{amsthm}
\usepackage{esint}
\usepackage{bm}
\usepackage{xcolor}
\usepackage{graphicx}
\usepackage[normalem]{ulem}
\usepackage{upgreek} 
\usepackage{mathtools}
\usepackage{mathrsfs}
\usepackage{subfiles}
\usepackage{enumitem}
\usepackage{hyperref}
\usepackage{accents}
\usepackage[normalem]{ulem}

\usepackage[mathlines]{lineno}

\usepackage{etoolbox}          %

\newcommand*\linenomathpatch[1]{%
  \cspreto{#1}{\linenomath}%
  \cspreto{#1*}{\linenomath}%
  \csappto{end#1}{\endlinenomath}%
  \csappto{end#1*}{\endlinenomath}%
}
\newcommand*\linenomathpatchAMS[1]{%
  \cspreto{#1}{\linenomathAMS}%
  \cspreto{#1*}{\linenomathAMS}%
  \csappto{end#1}{\endlinenomath}%
  \csappto{end#1*}{\endlinenomath}%
}

\expandafter\ifx\linenomath\linenomathWithnumbers
  \let\linenomathAMS\linenomathWithnumbers
  \patchcmd\linenomathAMS{\advance\postdisplaypenalty\linenopenalty}{}{}{}
\else
  \let\linenomathAMS\linenomathNonumbers
\fi

\linenomathpatch{equation}
\linenomathpatchAMS{gather}
\linenomathpatchAMS{multline}
\linenomathpatchAMS{align}
\linenomathpatchAMS{alignat}
\linenomathpatchAMS{flalign}

\makeatletter
\patchcmd{\mmeasure@}{\measuring@true}{
  \measuring@true
  \ifnum-\linenopenaltypar>\interdisplaylinepenalty
    \advance\interdisplaylinepenalty-\linenopenalty
  \fi
  }{}{}
\makeatother

\newcommand{\el}{\text{el}}
\newcommand{\vi}{\text{vi}}
\newcommand{\eq}{\text{eq}}
\newcommand{\ext}{\mathrm{ext}}

\newcommand{\mix}{\mathrm{mix}}

\newcommand{\R}{\mathbb{R}}

\newcommand{\N}{\mathbb{N}}

\newcommand{\M}{\mathbb{M}}
\newcommand{\calM}{\mathcal{M}}
\newcommand{\eps}{\varepsilon}
\newcommand{\pl}{\partial}
\newcommand{\DIV}{\mathrm{div\,}}

\newcommand{\bbM}{\mathbb{M}}
\newcommand{\Cof}{\mathrm{Cof}}

\newcommand{\wb}{\overline}
\newcommand{\wh}{\widehat}
\newcommand{\wt}{\widetilde}

\newcommand{\ubar}{\underline}
\newcommand{\dx}{\mathrm{\,d}x}
\newcommand{\dt}{\mathrm{\,d}t}

\newcommand{\dGamma}{\mathrm{\,d}S}%
\newcommand{\dd}{\mathrm{\,d}}
\newcommand{\dS}{\mathrm{\,d}S}
\newcommand{\id}{\mathrm{id}}

\newcommand{\Argmin}{\mathrm{Arg\,Min}}
\newcommand{\chibar}{\overline{\chi}}
\newcommand{\chihat}{\widehat{\chi}}
\newcommand{\mubar}{\overline{\mu}}

\newcommand{\cbar}{\overline{c}}
\newcommand{\chat}{\widehat{c}}

\newcommand\DT[1]{\mathchoice
                 {{\buildrel{\hspace*{.1em}\text{\LARGE.}}\over{#1}}}
                 {{\buildrel{\hspace*{.1em}\text{\Large.}}\over{#1}}}
                 {{\buildrel{\hspace*{.1em}\text{\large.}}\over{#1}}}
                 {{\buildrel{\hspace*{.1em}\text{\large.}}\over{#1}}}}
\newcommand\DDT[1]{\mathchoice
   {{\buildrel{\hspace*{.13em}\text{\LARGE.\hspace*{-.13em}.}}\over{#1}}}
   {{\buildrel{\hspace*{.1em}\text{\Large.\hspace*{-.1em}.}}\over{#1}}}
   {{\buildrel{\hspace*{.1em}\text{\large.\hspace*{-.1em}.}}\over{#1}}}
   {{\buildrel{\hspace*{.1em}\text{\large.\hspace*{-.1em}.}}\over{#1}}}}

\renewcommand{\dot}{\DT}
\renewcommand{\ddot}{\DDT}
\newcommand{\calE}{\mathcal{E}}

\newcommand{\calH}{\mathcal{H}}
\newcommand{\calR}{\mathcal{R}}

\newcommand{\rmD}{\mathrm{D}}

\newcommand{\mfh}{\mathfrak{h}}

\newcommand{\GL}{\mathrm{GL}}
\newcommand{\tr}{\mathrm{tr}}
\DeclarePairedDelimiterX{\abs}[1]{\lvert}{\rvert}{#1}
\DeclarePairedDelimiterX{\norm}[1]{\lVert}{\rVert}{#1}
\DeclarePairedDelimiterX{\ip}[1]{\langle}{\rangle}{#1}
\newcommand{\wto}{\xrightharpoonup{w}}
\newcommand{\wstarto}{\xrightharpoonup{w^*}}
\newcommand{\sto}{\xrightarrow{s}}

\theoremstyle{plain}
\newtheorem{Th}{Theorem}[section]
\newtheorem*{Th*}{Theorem}
\newtheorem{Lemma}[Th]{Lemma}
\newtheorem{Cor}[Th]{Corollary}
\newtheorem{Prop}[Th]{Proposition}

\theoremstyle{definition}
\newtheorem{Def}[Th]{Definition}
\newtheorem*{Def*}{Definition}

\newtheorem{Rem}[Th]{Remark}
\newtheorem{?}[Th]{Problem}
\newtheorem{Ex}[Th]{Example}

\newcommand{\tdots}{\,\vdots\,}

\title{Finite-strain poro-visco-elasticity\\ with degenerate mobility}
\author{Willem J. M. van Oosterhout\footnote{Weierstrass Institute for Applied Analysis and Stochastik, Mohrenstraße 39, Berlin, Germany, \url{willem.vanoosterhout@wias-berlin.de}, \url{matthias.liero@wias-berlin.de}} \and Matthias Liero$^\text{*}$}
\date{\today}

\numberwithin{equation}{section}
 \usepackage[width=16cm]{geometry}

\begin{document}

\maketitle
\begin{abstract}
    A quasistatic nonlinear model for poro-visco-elastic solids at finite strains is considered in the Lagrangian frame using the concept of second-order nonsimple materials. The elastic stresses satisfy static frame-indifference, while the viscous stresses satisfy dynamic frame-indifference. The mechanical equation is coupled to a diffusion equation for a solvent or fluid content. The latter is pulled-back to the reference configuration. To treat 
	the nonlinear dependence of the mobility tensor on the deformation gradient, the result by Healey \& Krömer is used to show that the determinant of the deformation gradient is bounded away from zero. Moreover, the focus is on the physically relevant case of degenerate mobilities. 
	 The existence of weak solutions is shown using a staggered time-incremental scheme and suitable energy-dissipation inequalities. 
\end{abstract}


\vspace{1em}

\section{Introduction}

The coupling of the  mechanical deformations of solids to other physical processes
such as heat conduction or diffusion of chemical species is relevant in many
applications in technology or biology.  We refer to \cite{MauBer1999MFFS,WLX*2020CPBT,Latr2004STPM,CaSeSa2022ECMS,HZZS2008TCDLDPG,CheAna2010CTFPLDEM} and the references therein for applications in thermo-mechanics, solid-state batteries, poroelasticity in biological tissue, hydrogen storage, and elastomeric materials. 

In many of these applications, the assumption of small strains is no longer justified and 
nonlinear, finite-strain, theories have to be considered. While the static theory for finite-strain elasticity 
developed rapidly after the seminal work 
of Ball \cite{Ball1977CIET}, 
the mathematical analysis for time-dependent processes in the case of large strains is just currently receiving increased attention,
see e.g.\ \cite{MiRoSa2018GERV}
for a result on finite-strain visco-plasticity, \cite{MieRou2020TVKVR,BaFrKr2023NaLMiTVE} 
for finite-strain thermo-visco-elasticity, and \cite{Roub2021CHEC,SchPes2023SILC} 
for phase-field models in the finite-strain setting.

In this text, we consider a coupled model for the visco-elastic evolution of the deformation $\chi(t,x)$ of a solid in the reference configuration $\Omega\subseteq\R^d$ and the diffusion of a solvent, fluid content, or other chemical species through the (visco-)elastic body. We shall call such materials 
\emph{poroelastic} (although this name is usually used to refer to the interaction between fluid flow and solids deformation, see \cite{Biot1941GTo3dC}).

 For a time horizon $T>0$ and $\Omega\subseteq \R^d$ a bounded, open reference configuration, we establish the existence of a deformation $\chi: [0,T]\times \Omega\to \R^d$ and a concentration $c:[0,T]\times \Omega\to \R^+$ satisfying the quasi-static system
\begin{subequations}
\label{Eqn:Intro}
\begin{alignat}{2} 
  \label{Eqn:Intro:Deformation}
- \DIV\!\big(\sigma_\text{el}(\nabla\chi,c)
+\sigma_\text{vi}(\nabla\chi,\nabla\dot\chi,c)-\DIV(\mathfrak{h}(\rmD^2\chi))\big) &= f(t),\quad &&\text{in}\ [0,T]\times\Omega,\\
   \label{Eqn:Intro:Diffustion}
\dot{c} - \DIV\!\big(\calM(\nabla\chi,c)\nabla \mu\big) &= 0,\quad &&\text{in}\ [0,T]\times\Omega,
\end{alignat}
\end{subequations}
where \emph{quasistatic} refers to the neglect of 
kinetic energy (and therefore  inertial forces $\rho \ddot{\chi}$) such that there 
will be no mechanical oscillations.
The total stress $\Sigma_\text{tot} := \sigma_\text{el}+\sigma_\text{vi}-\DIV \mathfrak{h}$ consists of the elastic stress $\sigma_\el(F,c) = \pl_{F} \Phi(F,c)$, 
coming from a free energy density $\Phi(\nabla\chi,c)$, the viscous stress $\sigma_\vi(F,\dot F,c) = \pl_{\dot F}\zeta(F,\dot F,c)$, 
given in terms of a viscous dissipation potential $\zeta(\nabla\chi,\nabla\dot \chi,c)$, and the hyperstress $\mathfrak{h}(G) = \partial_{G}\mathscr{H}(G)$ 
with potential $\mathscr{H}(\rmD^2\chi)$, which gives higher regularity of the deformation. 
Furthermore, $f$ is a body force, $\calM$ is the mobility tensor, 
and $\mu(F,c) = \pl_c\Phi(F,c)$ is the chemical potential. 
The system is completed by suitable boundary conditions. Note that the equation for the deformation is of fourth order, due to the 
hyperstress regularization. Thus, the Neumann boundary condition contains an additional
contribution related to the surface divergence of $\mathfrak{h}$, and a further boundary condition has to be satisfied, see  \eqref{Eqn:FiniteStrainBoundary2} and 
\eqref{Eqn:FiniteStrainBoundary3}. For the diffusion equation, we assume the following Robin-type boundary condition
\[
	\calM(\nabla\chi,c)\nabla\mu\cdot \Vec{n}  = \kappa(\mu_{\ext}(t){ -}\mu)\quad\text{on }\pl\Omega,
\]
where $\kappa\geq 0$ is a transmission coefficient and $\mu_\ext$ is a fixed external potential.
As is the case in previous works \cite{MieRou2020TVKVR}, \cite{RouTom2020DCEB}, we 
do not put any convexity assumptions on the energy density $\Phi$, but rather add the higher order convex term 
$\mathscr{H}(\rmD^2\chi)$, leading to the hyperstress $\mathfrak{h}$. 
This addition makes it into a second-grade non-simple material, a notion which was 
introduced by Toupin \cite{Toup1962EMwCS}.

We highlight that the diffusion processes takes place 
in the actual (deformed) configuration, which leads to 
a nontrivial description when pulled back to the reference 
configuration $\Omega$ for the mathematical analysis. 
In the actual configuration, the diffusion equation takes the 
form 
\begin{equation}
	\label{Eqn:Intro:EulerDiffusion}
\dot{ \mathsf{c}} + \mathrm{Div}(\mathsf{c} v - \M(F,\mathsf{c}) \nabla\upmu ) =0,
\end{equation}
where $\mathsf{c}(t,\chi(t,x)) = c(t,x)/\det \nabla\chi(t,x)$ is the spatial 
concentration, $v(t,\chi(t,x)) = \dot\chi(t,x)$ is the Eulerian velocity, $F(t,\chi(t,x)) = \nabla\chi(t,x)$ 
is the deformation gradient, $\upmu(t,\chi(t,x)) = \mu(t,x)$ is the spatial chemical potential,
and $\bbM$ is the Eulerian mobility tensor.
The Lagrangian mobility tensor $\calM$ is obtained from 
the pull-back of $\bbM$. The two are related
by the formula
\[
	\calM(F,c) = \frac{(\Cof F^\top) \M(F,c/\!\det F)\Cof F}{\det F}\quad \text{for }(F,c)\in \GL^+(d)\times \R^+,
\]
where $\Cof F = (\det F) F^{-\top}$ denotes the cofactor matrix associated with $F$.
Equation \eqref{Eqn:Intro:Diffustion} now follows from \eqref{Eqn:Intro:EulerDiffusion}
via the pullback under $\chi$.
As a consequence of the transformation, 
we need to ensure that the determinant of the deformation 
gradient is bounded away from zero, which follows from a result 
by Healey and Kr\"omer \cite{HeaKro2009IwSiSGNE}. Indeed, the 
hyperstress regularization in connection with  
growth assumptions for the elastic energy contribution in $\Phi$
in the form $\Phi(F,c) \geq C_1/(\det F)^q -C_2$, yields a 
uniform constant $\delta>0$ such that $\det \nabla\chi \geq \delta$
in $\Omega$ for all deformations with finite energy, see also \cite[Thm.~3.1]{MieRou2020TVKVR}. 
This blowup of the free energy, if the determinant of the deformation gradient approaches $0$ from above,
gives rise to local non-selfpenetration; however, we do not enforce global non-selfpenetration.

In \cite{DaNiSt2022ERMEM} a model for the growth of biological tissue was considered. Therein, the evolution of the 
deformation is coupled to the diffusion of a nutrient and the evolution of a growth variable which plays the role of a plastic tensor. However, the diffusion equation for the nutrient density is assumed to take place in the reference configuration and the deformation $\chi$ only appears in a source/sink term. In particular, 
no viscous stresses and no hyperstress regularization were considered. The missing temporal compactness of the deformation (gradient) is compensated by introducing a suitable time convoluted deformation $K*\chi$ instead of $\chi$ in the source terms. In \cite{RouSte2022*VSPS}, a similar system as \eqref{Eqn:Intro} was considered. However, there the authors 
work purely in the Eulerian setting, which restricts them to assuming that the actual domain is not evolving in time, i.e., $\chi(t,\Omega) = \Omega$.

Our new contribution is that we allow for degenerate mobilities, i.e.\ $\calM(F,0) \sim c^{m}$ for some power $m>0$. Until now, to the best of the authors' knowledge, the mobility was always assumed to be uniformly positive definite, see e.g.\ \cite{Roub2017VMSS,RouTom2020DCEB}, \cite{MieRou2020TVKVR} (which uses conductivities instead of mobilities), 
or \cite{Roub2021CHEC}. 
Our assumptions allow, for example, for the linear (actual) mobility $\M(c)=c\M_0$, or more generally, for a polynomial mobility. This type of mobility is physically relevant, as it models a higher species permeability whenever the material opens up (i.e., swells) due to an increase in species concentration, see \cite{CheAna2010CTFPLDEM}. See Example \ref{Ex:BiotFickDarcy} for the case of a Biot-type model.

Finally, an important modeling assumption is the static frame-indifference 
of the free energy, namely $\Phi(RF,c) = \Phi(F,c)$ for all $R\in\mathrm{SO}(d)$ and 
$F\in\GL^+(d)$, and the dynamic frame-indifference of the 
viscous stress potential, i.e., $\zeta(RF,\dot R F+R\dot F,c)=\zeta(F,\dot F,c)$ for all smooth $t\mapsto R(t)\in\mathrm{SO}(d)$,
and $F,\dot F\in \R^d$. 
For simplicity, we will assume 
that the viscous stress is linear in $(\nabla\chi^\top\nabla\chi)^{\dot{}}$, 
which is used to model non-activated dissipative processes with moderate rates, 
see e.g., \cite[Sect.\! 2]{MieRou2020TVKVR}. To control $\nabla\dot\chi$ via $(\nabla\chi^\top\nabla\chi)^{\dot{}}$, we exploit results for generalized Korn's inequalities by Neff \cite{Neff2002KFIN} and the extension by Pompe \cite{Pomp2003KFIV}. 
Here, again it is used that we can control the determinant of the deformation gradient via the hyperstress regularization.

The main result of this paper is that under certain conditions our model always admits weak solutions 
(see Definition \ref{Def:weakSolutionFiniteStrain} and Theorem~\ref{Thm:Main}). 
To prove the existence of solutions, we use a time-incremental scheme consisting of two steps. 
First, the elastic equation is solved to obtain the new deformation $\chi$, with fixed concentration $c$ from the previous time level, by reformulating the equation as minimization problem. 
Then, the diffusion equation is solved using a fixed point argument to obtain the new chemical potential $\mu$, which also implies the existence of a concentration $c$. Here, it is important to note that due to the degenerate mobility an extra regularization is needed, cf.\ \cite{Jung2015BbEMfCDS}. This regularization takes the form $\eta(-\Delta)^\theta \mu$, with regularization parameter $\eta>0$ and exponent $\theta > \frac{d}{2}$ (such that there is a compact embedding $H^\theta(\Omega)\hookrightarrow L^\infty(\Omega)$).
Next, we show that these solutions satisfy an energy-dissipation balance, from which suitable a priori estimates are derived. Due to the degenerate mobility, the classical Aubin--Lions lemma no longer suffices to pass to the limit in the concentration terms. Instead we use a more general compactness result due to Dubinski\u\i\ \cite{Dubi1965WCfNEaPE} to pass to the limit $(\tau,\eta)\to0$, see also Theorem \ref{Th:DubinskiiCompactness}.

The paper is structured as follows. In Section \ref{Sect:Setting} we introduce our model and assumptions in a mathematically rigorous way, and state the main result. In Section \ref{Sect:TimeDiscr}, we start the proof of the main result by introducing a regularization term and a time-incremental scheme. After showing existence of these approximate solutions, we proceed by showing an energy-dissipation inequality, from which suitable a priori estimates are deduced. Finally, in Section \ref{Sect:LimitPassage}, we pass to the simultaneous limit with respect to the time discretization and regularization to obtain the main result.

\section{Mathematical setting and main result}
\label{Sect:Setting}

Our model is described in the Lagrangian setting in the 
reference configuration $\Omega\subseteq \R^d$. We assume that 
$\Omega$ is an open, bounded domain with Lipschitz boundary, and that 
$\pl\Omega = \Gamma_D\cup\Gamma_N$ (disjoint) such that the Dirichlet part has positive surface measures $\int_{\Gamma_D} 1 \dS > 0$. We denote by $L^p(\Omega)$, $H^k(\Omega)$, and $W^{k,p}(\Omega)$ the usual Lebesgue and Sobolev spaces with the standard norms. Moreover, the (closed) subspace 
$W^{k,p}_0(\Omega)$ denotes the functions in $W^{k,p}(\Omega)$ with zero trace on $\Gamma_D$.
We consider deformations $\chi$ on $\Omega$ that are fixed on the Dirichlet part $\Gamma_D$, namely,
we consider the space
\[
W^{2,p}_\id(\Omega;\R^d) := \{\chi\in W^{2,p}(\Omega;\R^d) \mid \chi|_{\Gamma_D} = \id \}.
\] 
We denote by $``a\cdot b"$, $``\!A:B"$, and $``G \tdots H "$ 
the scalar products between vectors $a,b\in\R^d$, matrices $A,B\in \R^{d\times d}$, and third-order tensors $G,H\in\R^{d\times d\times d}$, respectively.

To introduce the system of partial differential equations coupling 
the evolution of the deformation $\chi$ and the concentration $c$ 
of some species, we consider a free energy density $\Phi=\Phi
(\nabla\chi,c)$, a higher-order regularization $\mathscr{H}
=\mathscr{H}(\rmD^2\chi)$, a viscous dissipation potential $\zeta 
= \zeta(\nabla\chi,\nabla\dot\chi,c)$, and a (Lagrangian) mobility 
tensor $\calM=\calM(\nabla\chi,c)$. The free energy density $\Phi$ 
gives rise to the first Piola--Kirchhoff stress $\sigma_\el$ 
and the chemical potential $\mu$, the viscous 
dissipation potential $\zeta$ to the viscous stress $\sigma_\vi$ via
\begin{equation}
\sigma_\text{el}(F,c) := \pl_F\Phi(F,c),\quad 
\mu(F,c) := \pl_c\Phi(F,c),\quad
\text{and}
\quad\sigma_\text{vi}(F,\dot F,c) := \pl_{\dot F}\zeta(F,\dot F, c),
\end{equation}
and the potential $\mathscr{H}$ 
to the hyperstress $\mfh(G) := \pl_G \mathscr{H}(G)$, where we have used 
the placeholders $F$ for $\nabla\chi$, $\dot F$ for $\nabla\dot 
\chi$, and $G$ for $\rmD^2\chi$. 
The model is then given by the evolutionary system in the reference domain $\Omega$
	\begin{subequations}
		\label{Eqn:FiniteStrainSystem}
		\begin{align} 
			\label{Eqn:FiniteStrainStress}
			- \DIV\!\big(\sigma_\el(\nabla\chi,c) + \sigma_\vi(\nabla\chi,\nabla\dot\chi,c) - \DIV \mfh(\rmD^2 \chi)\big) &= f(t),\\
			\label{Eqn:FiniteStrainChemPot}
			\dot{c} - \DIV\!\big(\calM(\nabla\chi,c)\nabla \mu\big) &= 0,
		\end{align} 
	\end{subequations}
	completed with the boundary conditions
	\begin{subequations} 
		\label{Eqn:FiniteStrainBoundary}
		\begin{alignat}{2}
			\label{Eqn:FiniteStrainBoundary1}&\chi = \id &\text{on}\ \Gamma_D,\\
			\label{Eqn:FiniteStrainBoundary2}&\big(\sigma_\el(\nabla\chi,c) + \sigma_\vi(\nabla\chi,\nabla\dot\chi,c)\big)\Vec{n} - \DIV_{\!\text{s}}(\mfh(\rmD^2\chi)\Vec{n}) = g(t) \hspace{25pt}&\text{on}\ \Gamma_N,\\
			\label{Eqn:FiniteStrainBoundary3}&\mfh(\rmD^2\chi):(\Vec{n}\otimes\Vec{n}) = 0 &\text{on}\ \pl\Omega,\\
			\label{Eqn:FiniteStrainBoundary4}&\calM(\nabla\chi,c)\nabla\mu\cdot \Vec{n}  = \kappa(\mu_{\ext}(t){ -}\mu) &\text{on}\ \pl\Omega,
		\end{alignat}
	\end{subequations} 
	where $\vec{n}$ denotes the unit normal vector on $\pl\Omega$, and
	$\kappa(x)\geq 0$ and $\mu_{\ext}(t,x)$ are a given permeability and an external potential, respectively. 	Here, $\DIV\!_\text{s}$ denotes the surface divergence, defined by $\DIV\!_\text{s}(\cdot) = \tr(\nabla_\text{s}(\cdot))$, i.e., the trace of the surface gradient $\nabla_\text{s} v = (I-\Vec{n}\otimes\Vec{n})\nabla v = \nabla v - \frac{\pl v}{\pl\Vec{n}}\Vec{n}$. Finally, we consider initial conditions
	\begin{equation}
		\label{Eqn:FiniteStrainInitial}
		\chi(0) = \chi_0,\quad c(0) = c_0 \quad\text{on}\ \Omega.
	\end{equation}

For any $R>0$, let us denote the set 
\begin{equation}\label{eq:defSetFR}
	\mathsf{F}_R:= \big\{F\in \GL^+(d)\,\big|\,|F|\leq R,~|F^{-1}|\leq R,\text{ and }\det F \geq 1/R\big\}.
\end{equation}

	We impose the following assumptions on our model:

\begin{enumerate}[label=\emph{\bfseries (A\arabic*)}]
\item \label{Assu:Hyperstress} The hyperstress potential is a convex, frame-indifferent $C^1$ 
function $\mathscr{H}:\R^{d\times d\times d}\to \R^+$ such 
that the hyperstress is given by $\mfh(G) = \partial_G \mathscr{H}(G)\in\R^{d\times d\times d}$. Moreover, %
there exist $p\in (d,\infty)\cap [3,\infty)$ and constants $C_{\calH,1}$, $C_{\calH,2}$, $C_{\calH,3}>0$ such that 
\[
C_{\calH,1}\abs{G}^p\leq \mathscr{H}(G)\leq C_{\calH,2}(1+\abs{G}^p), \quad\abs{\pl_G\mathscr{H}(G)}\leq C_{\calH,3}\abs{G}^{p-1}\quad \text{for all}\  G\in \R^{d\times d\times d}.
\]

\item \label{Assu:MobilityTensor} 
The mobility tensor $\calM:\GL^+(d)\times \R^+
\to\R^{d\times d}_{\text{sym}}$ is a continuous map.
There exist an exponent $m>0$, and for all $R>0$ there exist constants $C_{0,\calM,R},C_{1,\calM,R}>0$ such that
\begin{equation}
	\begin{gathered}
\xi\cdot \calM(F,c)\xi\geq C_{0,\calM,R} c^{m}\abs{\xi}^2 
\quad\text{and}\quad 		|\calM(F,c)|\leq C_{1,\calM,R} c^{m}\\
\text{for all}\ \xi\in\R^d, ~F\in \mathsf{F}_R, ~c\in \R^+.
	\end{gathered}
\end{equation}
The admissible range of the exponent $m>0$ depends on the 
growth properties of (the derivatives of) $\Phi$ and is fixed in assumption 
\ref{Assu:FreeEn}.

\item \label{Assu:FreeEn} The free energy $\Phi: \GL^+(d)\times\R^+ \to \R$ is bounded from below, continuous, and $C^2$ on $\GL^+(d)\times (0,\infty)$, i.e., for strictly positive concentrations. It is frame indifferent, and satisfies the following assumptions:
\begin{itemize}
    \item[(i)] For any $c\in \R^+$ there exists constants $C_{\Phi,0}, C_{\Phi,1}>0$ such that
    \begin{equation*}
\Phi(F,c) \geq C_{\Phi,0}\abs{F} + \frac{C_{\Phi,0}}{(\det F)^q}-C_{\Phi,1}\quad \text{for all}\ F\in \GL^+(d).
\end{equation*}
	\item[(ii)] There exist an exponent $-1<r<\infty$ such that $r+m\geq 0$, and for all $R>0$ constants $C_i: = C_{\Phi,i,R}>0$ ($1\leq i\leq 2$) and constants $\gamma_i := \gamma_{\Phi,i,R}\geq 0$ ($1\leq i\leq 2$) such that 
	\[
	\frac{C_1}{c} + \gamma_1c^r \leq \pl^2_{cc}\Phi(F,c) \leq \frac{C_3}{c} + \gamma_2c^r \quad \text{for all}\ c\in \R^+,~F\in \mathsf{F}_R.
	\]
	Concerning the constants $\gamma_i$, we distinguish two cases (see also Remark \ref{Rem:DiffCaseI/II}): 
 \begin{description}
        \item{\textbf{Case I}:}  We assume $\gamma_1=\gamma_2 = 0$, and also require $1 \leq m\leq 2$. 
        \item{\textbf{Case IIa}:} For $\gamma_2\geq \gamma_1>0$, we require $0<m\leq 3+r$.
        \item{\textbf{Case IIb}:} For $\gamma_2\geq \gamma_1>0$, we require $0<m\leq 2$.
\end{description}

	\item[(iii)] There exist an exponent $\alpha\in \R$, and for all $R>0$ a constant $C_{\Phi,5,R}>0$ such that 
	\[
	\abs*{\pl^2_{Fc}\Phi(F,c)}\leq C_{\Phi,5,R}c^\alpha \quad \text{for all}\ c\in \R^+,~ F\in\mathsf{F}_R
	\]
	   In \textbf{Case I} above, 
        $\alpha$ is such that $0\leq m+\alpha \leq \frac{p-s}{ps}$, where $1<s=\frac{md+2}{md+1}< 2$ and $0\leq m+2\alpha$.

    In both \textbf{Case IIa} and \textbf{Case IIb} $\alpha\geq -1$ is such that $0\leq m+\alpha \leq (2+r)\frac{p-s}{ps}$, where $1<s=\min\{\frac{md+2(r+2)}{md+r+2},\frac{d(m+r+1)+2(r+2)}{d(m+r+1)+r+2}\}<2$. Furthermore, in case \textbf{Case IIa} we require that $0\leq m+2\alpha <m+1+r$, while in \textbf{Case IIb} we require that $0\leq m+2\alpha <m+2+2r$

\end{itemize}

\item \label{Assu:EnergyInitFin} For all $R>0$ there exists a concentration $c_R\in \R^+$ such that $\Phi(F,c_R)<\infty$ and $|\pl_c\Phi(F,c_R)|<\infty$ for all $F\in \mathsf{F}_R$.		
	
\item \label{Assu:ViscousStress} The viscous stress potential 
$\zeta:\R^{d\times d}\times \R^{d\times d}\times \R^+\to \R^+$ 
is such that 
$\zeta(F,\dot F,c) = \wh\zeta(F^\top F,(F^\top F)^{\dot{}} ,c)$, where $\wh\zeta:\R^{d\times d}_\text{sym} \times \R^{d\times d}_\text{sym}\times \R^+ \to \R^+$
is quadratic in the second variable, namely
\[
\wh\zeta(\mathsf{C},\dot{\mathsf{C}},c) 
= \frac{1}{2}\dot {\mathsf{C}}:\mathbb{\wt D}(\mathsf{C},c)\dot {\mathsf{C}}.
\]
This quadratic form is such that 
there exist constants $C_{\zeta,1}$, $C_{\zeta,2}>0$ such that
\[
C_{\zeta,1}\abs{(F^\top F)^{\dot{}} }^2 \leq \wh\zeta(F^\top F,(F^\top F)^{\dot{}},c) \leq C_{\zeta,2}\abs{(F^\top F)^{\dot{}} }^2 \quad \text{for all}\ c\in \R^+,~F\in \R^{d\times d}.
\]

\item \label{Assu:Loading}  The external forces satisfy $f\in W^{1,\infty}(0,T;L^{2}(\Omega;\R^d))$, $g\in W^{1,\infty}(0,T;L^{2}(\pl\Omega;\R^d))$. We set
\[
\ip{\ell(t),\chi} := \int_\Omega f(t)\cdot \chi\dx + 
\int_{\Gamma_N} g(t)\cdot\chi \dS
\] 
such that $\ell\in W^{1,\infty}(0,T;H^{1}(\Omega;\R^d)^*)$.

\item \label{Assu:Permeability} The permeability $\kappa \in L^\infty(\partial\Omega)$ is nonnegative and strictly positive on a part of the boundary $\pl\Omega$ with positive surface measure, i.e., $\int_{\pl\Omega}\kappa\dS \geq \kappa_* > 0$. We assume that the external chemical potential is such that $\mu_\ext\in L^2([0,T]{\times}\partial\Omega)$.

\item \label{Assu:RegularityInit} The initial conditions satisfy $\chi_0\in W^{2,p}_\id(\Omega;\R^d)$ with $\det\nabla\chi_0\geq \rho_0>0$ and $c_0\in L^1(\Omega)$ with $c_0\geq 0$ and are such that $\int_\Omega \Phi(\nabla\chi_0,c_0)\dd x< \infty$.
\end{enumerate}
	
As an immediate consequence of these assumptions, we note the following lower bound for the free energy.			
\begin{Lemma}[Lower bound for free energy]
\label{Lemma:BelowEstimatePhiMix}
If %
$\nabla\chi\in \mathsf{F}_R$ for some $R>0$, then there exist $C_1, C_2, C_3>0$ such that
\begin{equation}
	\label{Eqn:LowerBoundPhiMix}
\int_\Omega \Phi(\nabla\chi,c)\dx \geq C_1\norm{c}_{L\log L(\Omega)} + \gamma_1C_2\norm{c}_{L^{2+r}}^{2+r} - C_3.
\end{equation}
\end{Lemma}		
\begin{proof}
Using first the lower bound $\pl_{cc}^2\Phi(F,c)\geq \frac{C_{\Phi,1,R}}{c}$ of Assumption \ref{Assu:FreeEn}(ii) and Assumption \ref{Assu:EnergyInitFin}, we obtain by integrating once 
\[
\pl_c\Phi(F,c) \geq C_{1,\Phi,R}\log(c) + \pl_c\Phi_\mix(F,c_R) - C_{1,\Phi,R}\log(c_R).
\]
Integrating again, and using Young's inequality with $\epsilon$ to absorb the lower order terms in $c$ in the highest order term, we thus obtain
\[
\int_\Omega \Phi(\nabla\chi,c)\dx \geq C_1\norm{c}_{L\log L(\Omega)} - C_3.
\]
The second part of the lower bound now follows from integrating the lower bound $\pl^2_{cc}\Phi \geq \gamma_1c^r$ twice.
\end{proof}

    \begin{Rem}%
		\begin{enumerate}[label=(\roman*),ref=\theRem(\roman*)]
    \item(Growth of free energy) Note that Assumption \ref{Assu:FreeEn}(ii) implies that $c\mapsto \Phi(c,F)$ is superlinear for every $F\in\GL^+(d)$. Moreover, we see that $\Phi(F,c) \sim c\log c -c +\gamma_1 c^{r+2}$.

\item\label{Rem:DiffCaseI/II}(Differences between Case I and II) We remark that the limit $r \downarrow-1$ in the conditions of Case II almost reduces to the conditions of Case I. However, 
in Case I, we additionally need to assume the lower bound $m\geq 1$. Without this bound, it is not possible to obtain the strong convergence of the concentrations, which is needed to pass to the limit, see Prop.~\ref{Prop:LimitTimeDisc}. Note that as a result, Case I is restricted to negative exponents $\alpha< 0$. Indeed, we have $\alpha\leq \frac{p(1-ms)-s}{ps} <0$
since $m\geq 1$ and $s>1$. In fact,
for $p\to \infty$, we get $\alpha\leq \frac{1}{s}-m$.

    \item\label{Rem:RelationMuC}(Relation between $\mu$ and $c$) From Assumption \ref{Assu:FreeEn}(ii) and the resulting strict monotonicity of $c\mapsto \pl_c\Phi(F,c)$ 
	it follows that the equation $\mu = \pl_c\Phi(F,c)$ for $F\in\GL^+(d)$ fixed is uniquely solvable for given $\mu\in\R$. By the implicit function theorem, the map $c=c(F,\mu)$ is continuous on $\GL^+(d)\times \R$. In particular, for any $F\in\GL^+(d)$ and $\mu\in\R$ there exist constants $M,\delta>0$ depending only on 
	$|\mu|$, $|F|$, and $\det F$ such that $\delta \leq c(F,\mu)\leq M$. 
    
    \item(Constraint $c\geq 0$)    To obtain the constraint that the concentration stays positive, we note that $\pl^2_{cc} \Phi(F,c)\geq \frac{C_1}{c}$ (cf. \ref{Assu:FreeEn}(ii)) implies that $\lim_{c\to 0^+}\pl_c \Phi(F,c) = -\infty$. Thus, we can extend $\Phi$ by setting $\Phi(F,c) = +\infty$ for $c\leq 0$. Note that this implies that $\pl_c\Phi(F,c) = \emptyset$ for $c\leq 0$ (and in particular for $c=0$).
    
    Another possibility (which we will not use) is discussed in \cite[Rem.~2]{Roub2017ECTDfPEM}, where instead the mobility is extended by setting $\calM(F,c)=0$ for $c<0$. However, this will only give the positivity of the concentration for solutions \eqref{Eqn:FiniteStrainSystem}, not for the regularized problem. Indeed, since our existence proof uses an extra regularization for the diffusion part, it is not clear if the regularized time-discrete solution $c_k$ obtained from \eqref{Eqn:TimeDisc:FiniteStrainSystem} is positive. Consequently, several modifications would have to be made to the proof.
    
    \item(Degenerate \& non-degenerate mobility)
		If we assume that the mobility is uniformly positive definite, continuous and bounded, i.e., if $m=0$ (see \cite{Roub2017VMSS}, \cite{MieRou2020TVKVR} (conductivities instead of mobilities), or \cite{RouTom2020DCEB}), then it is possible to prove existence of solutions without using an additional regularization. Indeed, in this case one can test the diffusion equation with $\mu$ and use the estimate $\calM \nabla\mu\cdot\nabla\mu \geq C\abs{\nabla\mu}^2$ to obtain a bound for $\nabla\mu$ in $L^2(\Omega)$. Here, we restrict the discussion to the degenerate case, where $m>0$. %

    \item(Gradient structure) The viscous-elastic evolution and the diffusion process can be formally written
in terms of a gradient-flow structure. Indeed, define the free energy functional
$\calE_0(\chi,c) = \int_{\Omega}\Phi(\nabla\chi,c)+\mathscr{H}(\rmD^2\chi)\dd x$
as well as the convex dissipation potential
\begin{multline*}
\calR_\text{tot}(\chi,c,\dot\chi,\dot c)
= \int_\Omega \zeta(\nabla\chi,\nabla\dot \chi,c)\dd x \\+ \int_\Omega \calM(\nabla\chi,c)\nabla({-}\Delta_{\calM(\nabla\chi,c),\kappa}^{-1}\dot c)\cdot \nabla({-}\Delta_{\calM(\nabla\chi,c),\kappa}^{-1}\dot c)\dd x + \int_\Omega \kappa ({-}\Delta_{\calM(\nabla\chi,c),\kappa}^{-1}\dot c)^2 \dd\mathrm{S}
\end{multline*}
with $-\Delta_{\calM(\nabla\chi,c),\kappa}^{-1} : v \mapsto \tilde\mu$ denoting the linear operator defined formally by the weak solution $\tilde\mu$ of the elliptic equation $-\DIV(\calM(\nabla\chi,c)\nabla \tilde\mu)=v$
with boundary conditions $\calM(\nabla\chi,c)\nabla\mu\cdot\Vec{n} + \kappa\tilde\mu=0$. The evolutionary system in \eqref{Eqn:Intro} can be formally rewritten in the form
\[
\pl_{(\dot\chi,\dot c)}\calR_\text{tot}(\chi,c,\dot\chi,\dot c) + \rmD_{(\chi,c)}\calE_0(\chi,c) =  \xi_\text{ext},	
\]
where $\xi_\text{ext}$ contains the mechanical loading and the
external potential $\mu_\text{ext}$. It is possible to exploit the gradient 
structure of the equations to obtain time-discrete solutions via a staggered incremental scheme. However, since this scheme would require additional assumptions to derive suitable a priori estimates, we instead use a fixed point theorem to obtain the existence of time-discrete concentrations.
\end{enumerate}
\end{Rem}

    \begin{Ex}[Strong coupling in free energy] Consider the multiplicative decomposition of the deformation gradient as $F=F_\el F_c$, see e.g. \cite{Luba2004CTMD}. In this case, the free energy is formulated in terms of $F_\el = FF_c^{-1}$, e.g., $\Phi(F,c) = \Phi_1(FF^{-1}_c)+ \Phi_2(c)$, where $\Phi_1:\GL^+(d)\to [0,\infty)$ and $\Phi_2:(0,\infty)\to [0,\infty)$ describe the elastic and the chemical contribution to the free energy, respectively.
    A standard choice for $F_c$, which models isotropic swelling, 
    is e.g., $F_c = a(c)I$ for some $C^1$ function $a:[0,\infty)\to(0,\infty)$ with $0<a_*\leq a(c)\leq a^*$.
    
    We assume that $\Phi_1(F_\el,c)$ satisfies the coercivity estimate $\Phi_1(F_\el,c)\geq C_0\abs{F_\el} + \frac{C_0}{(\det F_\el)^q} - C_1$ with exponent $1<q<\infty$ as in Assumption~\ref{Assu:FreeEn}(i).
    Then, it follows that 
    \[
    \Phi_1(\nabla\chi F_c^{-1},c) 
\geq \frac{C_0}{a(c)}\abs{\nabla\chi} + \frac{C_0a(c)^q}{(\det F_\el)^q} - C_1\geq \frac{C_0}{a^*}\abs{F_\el} + \frac{C_0a_*}{(\det F_\el)^q} - C_1, 
    \]
    i.e., $\Phi(F,c)$ also satisfies the coercivity estimate. However, conditions on $\Phi_1$ and $\Phi_2$ for convexity with respect to $c$ and the growth conditions in Assumption \ref{Assu:FreeEn}(ii) and (iii) are hard to formulate at this abstract level and have to be confirmed for more concrete examples.

    \end{Ex}

	\begin{Ex}[Biot model and Fick/Darcy's law]
	\label{Ex:BiotFickDarcy}
	The Biot model \cite{Biot1941GTo3dC} (cf.\ \cite[Sect.~4]{RouTom2020DCEB})
	with Boltzmann entropy is given as	
	\begin{equation*}
	\Phi(F,c) = \Phi_\el(F) + \frac{1}{2}M_B(c-c_\eq-\beta(\det F{-}1))^2 
	+ k c\big(\log\big(\frac{c}{c_\eq}\big)-1\big),
	\end{equation*}
	for some suitable elastic energy $\Phi_\el$ and constants $M_B,\beta,k, c_\eq>0$. In this case, the assumptions for Case II are satisfied with $\alpha=0$ and $r=0$. Defining the (Eulerian) flux as $\boldsymbol{j}=-\M(F,c)\nabla\mu$, and assuming the mobility is linear in $c$, namely $\M(F,c)=c\M_0$ (i.e., $m=1$), we then obtain
	\[
		\boldsymbol{j}=-k\M_0 \nabla c - c\M_0\nabla p, 
	\]
	where $p=M_B(c-c_\eq - \beta(\det F{-}1))$ is the pressure. The first term corresponds to Fick's law, while the second is related 
	to Darcy's law.
\end{Ex}

	\begin{Def}[Weak solution]
		\label{Def:weakSolutionFiniteStrain}
		Let $1<s<2$ be as in \ref{Assu:FreeEn}(iii). We call a pair $(\chi,c)$ a weak solution of the  
		initial-boundary-value problem \eqref{Eqn:FiniteStrainSystem}--\eqref{Eqn:FiniteStrainInitial}
		if $\chi\in L^\infty(0,T;W^{2,p}_\id(\Omega;\R^d))$, $\dot\chi\in L^2(0,T;H^1(\Omega;\R^d))$ and $c\in L^\infty(0,T;L\log L(\Omega))$, $\dot c \in L^s(0,T;W^{1,s}(\Omega)^*)$ with $\nabla c^{\frac{m}{2}}\in L^2(0,T;L^2(\Omega))$ (Case I). In Case II, we additionally require that $c\in L^\infty(0,T;L^{2+r}(\Omega))$ and $\nabla c^{\frac{m+1+r}{2}}$, $\nabla c^{\frac{m}{2}+1+r}\in L^2(0,T;L^2(\Omega))$. The pair satisfies the integral equations 
		\begin{subequations}
        \begin{equation}
		\label{Eqn:FiniteStrain:WeakSoln:chi}
		\int_0^T\int_\Omega \big(\sigma_\el(\nabla\chi,c) + \sigma_\vi(\nabla\chi,\nabla\dot\chi,c)\big) : \nabla \phi 
		+ \mfh(\rmD^2\chi)\tdots \rmD^2\phi\dx\dt 
		= \int_0^T \ip{\ell(t),\phi}\dt
		\end{equation} 
		\end{subequations}
        for all %
		$\phi\in L^2(0,T;W_0^{2,p}(\Omega;\R^d))$, where $\ip{\cdot,\cdot}$ denotes the duality pairing between $W^{2,p}(\Omega;\R^d)^*$ and $W^{2,p}(\Omega;\R^d)$, %
        and
		\begin{equation}
		\label{Eqn:FiniteStrain:WeakSoln:mu}
		\int_0^T \ip{\dot c,\psi}\dt + \int_0^T\int_\Omega\calM(\nabla\chi,c)\nabla\mu \cdot \nabla \psi \dx\dt + \int_0^T\int_{\pl\Omega} \kappa(\mu - \mu_\ext)\psi\dGamma\dt = 0
		\end{equation} 
		for all $\psi\in L^{s'}(0,T;W^{1,s'}(\Omega))$, where $\ip{\cdot,\cdot}$ denotes the duality pairing between $W^{1,s}(\Omega)^*$ and $W^{1,s'}(\Omega)$.
  
		Furthermore, we require that $\mu\in \pl_c\Phi(\nabla\chi,c)$ almost everywhere in $\Omega$, and that $\mu\in L^2([0,T]\times\pl\Omega)$.
	\end{Def}
	
    We note that sufficiently smooth weak solutions indeed satisfy the classical formulation \eqref{Eqn:FiniteStrainSystem} with boundary conditions \eqref{Eqn:FiniteStrainBoundary}. For the derivation of \eqref{Eqn:FiniteStrainStress} and \eqref{Eqn:FiniteStrainBoundary1}-\eqref{Eqn:FiniteStrainBoundary3}, we refer to \cite[Eqn. (2.28)-(2.29)]{MieRou2020TVKVR}. The derivation of \eqref{Eqn:FiniteStrainChemPot} and \eqref{Eqn:FiniteStrainBoundary4} follows directly from integration by parts.
 
	\begin{Rem}
 Note that we not state any regularity conditions for the chemical potential $\mu$, but only for the concentration $c$. 
 Using the relation $\nabla\mu = \pl^2_{Fc}\Phi\rmD^2\chi + \pl^2_{cc}\Phi\nabla c$ and the bounds in \ref{Assu:MobilityTensor} and \ref{Assu:FreeEn}, we see that this gives a well-defined concept of weak solution.%
	\end{Rem}	
	
	\begin{Th}[Existence of weak solutions]
		\label{Thm:Main}
	Suppose that the assumptions \ref{Assu:Hyperstress}--\ref{Assu:RegularityInit} hold. Then, the system in \eqref{Eqn:FiniteStrainSystem}--\eqref{Eqn:FiniteStrainInitial} possesses at least 
	one weak solution in the sense of Definition \ref{Def:weakSolutionFiniteStrain}.
	\end{Th}

\section{Time-discretization of a regularized problem}
\label{Sect:TimeDiscr}

	For the time discretization, we consider an equidistant partition $\{0=t_0<t_1<...<t_{N}=T\}$ of $[0,T]$, where 
	$N\in\mathbb{N}$, $t_k = k\tau$ with $k=0,\ldots,N$ and $\tau=T/N>0$. 
	We construct approximate solutions $(\chi_k,\mu_k)$ such that $\chi_k \approx \chi(t_k)$ and $\mu_k\approx \mu(t_k)$ using a staggered scheme.
In the following, we use for the time discretization the difference notation
\[
\delta_\tau h_k = \frac{h_k-h_{k-1}}{\tau} \quad\text{for}\quad k =1,\ldots,N.
\]
We consider the following time-discrete system for $k=1,\ldots, N$:
\begin{subequations}
    \label{Eqn:TimeDisc:FiniteStrainSystem}
    \begin{align} 
			\label{Eqn:TimeDisc:FiniteStrainStress}
			\int_\Omega \big(\sigma_\el(\nabla\chi_k,c_{k-1}) + \sigma_\vi(\nabla\chi_{k-1},\delta_\tau\nabla\chi_k,c_{k-1})\big) : \nabla\phi 
			+ \mfh(\rmD^2 \chi_k) \tdots \rmD^2\phi\dx &= \ip{\ell_k,\phi},
			\\
			\label{Eqn:TimeDisc:FiniteStrainChemPot}
			\ip{\delta_\tau c_k,\psi} + \int_\Omega  \calM(\nabla\chi_{k},c_{k})\nabla\mu_k \cdot \nabla\psi \dx%
			+\int_{\pl\Omega}\kappa(\mu_k{-}\mu_{k,\ext})\psi\dS &= 0,\\[0.35em]
			\text{where }\mu_k \in \partial_{c}\Phi(\chi_k,c_k)\text{ a.e.\ in }\Omega.&
		\end{align} 
   Moreover, we set $\ell_k := \frac{1}{\tau}\int_{(k-1)\tau}^{k\tau}\ell(t)\dt$ and  $\mu_{k,\ext}:= \frac{1}{\tau}\int_{(k-1)\tau}^{k\tau}\mu_\ext(t)\dt$ for $k=1,\ldots N$ and set $\ell_0:=\ell_1$ such that $\delta_\tau \ell_1 =0$.

We note that this scheme is constructed in a specific way: Starting from the initial 
conditions $\chi_{k=0} = \chi_0$ and $c_{k=0} = c_0$ (cf.~\ref{Assu:RegularityInit}), 
we first solve the mechanical equation using the concentration $c_{k-1}$ from the previous time step, 
and then we solve the diffusion equation using the just obtained deformation $\chi_{k}$. 
The existence of time-discrete solutions $(\chi_k,c_k)$ follows from the formulation of the mechanical equation as
the Euler--Lagrange equation of a minimization problem 
and for the diffusion equation via a fixed-point argument.

Due to the degeneracy of the mobility, we cannot expect to obtain the existence and
suitable estimates for $\nabla\mu_k$. Thus, following \cite{Jung2015BbEMfCDS}, we add an elliptic 
regularization $\eta(-\Delta)^\theta\mu$ to the diffusion equation, 
for some small parameter $0<\eta\ll 1$ and $\theta\in\N$. 
This regularization ensures that $\mu_k\in H^\theta(\Omega)$. 
The exponent $\theta$ is chosen sufficiently large such that 
$H^\theta(\Omega)\hookrightarrow L^{\infty}(\Omega)$ compactly, i.e., $\theta>\frac{d}{2}$. 
We will use the notation $|\beta|:=\sum_{i=1}^d \beta_i$ for multi-indices $\beta\in\N_{0}$.
With this regularization, we replace the equation in \eqref{Eqn:TimeDisc:FiniteStrainChemPot} with 
\addtocounter{equation}{-1}
\begin{multline}
\tag{\theequation*}
\label{Eqn:TimeDisc:FiniteStrainChemPot:Reg}
\ip{\delta_\tau c_k,\psi} + \int_\Omega  \calM(\nabla\chi_{k},c_{k})\nabla\mu_k \cdot \nabla\psi \dx \\+ \eta\int_\Omega\sum_{\abs{\beta}=\theta}\pl^\beta\mu_k\cdot \pl^\beta\psi \dx + \int_{\pl\Omega}\kappa(\mu_k{-}\mu_{k,\ext})\psi\dS = 0
\end{multline}
	\end{subequations}
for all $\psi\in H^\theta(\Omega)$ and with $\mu_k \in \partial_{c}\Phi(\chi_k,c_k)$ a.e.\ in $\Omega$.
Note that this weak formulation is still well-defined since 
$c_k\in L^\infty(\Omega)$ (see Remark~\ref{Rem:RelationMuC}) and
thus $\calM(\nabla\chi_{k},c_{k})\in L^\infty(\Omega)$ and $\nabla\mu_k$, $\nabla\psi\in L^2(\Omega)$.

\begin{Rem}
We highlight that we do not require an
additional regularization term $\eta \delta_\tau \chi_k$ in the mechanical equation
to deal with the frame-indifferent viscous dissipation potential $\zeta$ as e.g.\ in \cite{MieRou2020TVKVR}.
\end{Rem}

\begin{Lemma}[Existence of time-discrete solutions]
\label{lemma:DiscreteSolutionsExistence}
Let the assumptions in \ref{Assu:Hyperstress}--\ref{Assu:RegularityInit} hold. We fix $\mu_0^\eta\in H^\theta(\Omega)$
and $\chi_0^\eta\in W^{2,p}_\id(\Omega;\R^d)$ such that $\det\nabla\chi_0^\eta\geq \rho_0>0$, and 
we define $c_0^\eta$ as the unique density with $\mu_0^\eta =\pl_c\Phi(\nabla\chi_0^\eta,c_0^\eta)$. 
Then, starting for $k=0$ from $\chi_0^\eta$ and $c_0^\eta$, we can iteratively find weak solutions $(\chi_k^\eta,c_k^\eta,\mu_k^\eta)\in W^{2,p}_\id(\Omega;\R^d)\times L^{\infty}(\Omega)\times H^\theta(\Omega)$, 
 for $k=1,\ldots, N$, 
for the mechanical equation \eqref{Eqn:TimeDisc:FiniteStrainStress} 
and the regularized diffusion equation
\eqref{Eqn:TimeDisc:FiniteStrainChemPot:Reg}, respectively,   
with $c_k^\eta> 0$
and $\mu_k^\eta=\pl_c\Phi(\nabla\chi_k^\eta,c_k^\eta)$ a.e.\ in $\Omega$.
\end{Lemma}

\begin{proof} We drop the index $\eta$ throughout the proof.

\textit{Mechanical step.} First, we consider the equation  \eqref{Eqn:TimeDisc:FiniteStrainStress}, where $\chi_{k-1}\in W^{2,p}_\id(\Omega;\R^d)$ and $c_{k-1}\in L^{\infty}(\Omega)$ with $c_{k-1}>0$ 
a.e.\ in $\Omega$ are fixed. 
Let us define the energy functional $\calE_0$ and the dissipation potential $\calR$ via
\[
\calE_0(\chi,c) := \int_\Omega \Phi(\nabla\chi,c)+\mathscr{H}(\rmD^2\chi)\dx,
\quad\text{and}\quad
\calR(\chi,\dot\chi,c) := \int_\Omega \zeta(\nabla\chi,\nabla\dot\chi,c)\dx.
\]
We assume that $\chi_{k-1}$ and $c_{k-1}$ are such that $\calE_0(\chi_{k-1},c_{k-1})<\infty$.
The deformation $\chi_k\in W^{2,p}_\id(\Omega;\R^d)$ is obtained
as the solution of the minimization problem 
\begin{equation}
\label{Eqn:TimeDisc:MinimizerChi}
\chi_k \in \Argmin\Big\{ \mathcal{E}_0(\wt\chi,c_{k-1}) + \tau\calR\Big(\chi_{k-1},\frac{\wt\chi-\chi_{k-1}}{\tau},c_{k-1}\Big) - \ip*{\ell_k,\wt\chi} \,\Big|\,\wt{\chi}\in W^{2,p}_\id(\Omega;\R^d)\Big\}.
\end{equation} 
Note that the corresponding weak Euler--Lagrange equation 
is exactly \eqref{Eqn:TimeDisc:FiniteStrainStress}. Indeed, the G\^ateaux differentiability
of $\chi \mapsto \calE_0(\chi,c_{k-1})$ follows as in \cite[Proposition 3.2]{MieRou2020TVKVR}.

To show that there exists a minimizer $\chi_k\in W^{2,p}_\id(\Omega;\R^d)$, we first note that the functional in \eqref{Eqn:TimeDisc:MinimizerChi} is bounded from below due to the assumptions \ref{Assu:Hyperstress}, and \ref{Assu:FreeEn}--\ref{Assu:Loading}. This boundedness follows from the embedding $W^{2,p}_\id(\Omega)\hookrightarrow W^{1,\infty}(\Omega;\R^d)$, which gives a bound for $\nabla\chi$ in $L^\infty(\Omega;\R^{d\times d})$. %
In fact, using the assumptions on the hyperstress, we also obtain that the functional in \eqref{Eqn:TimeDisc:MinimizerChi} is coercive on $W^{2,p}_\id(\Omega;\R^d)$. To show that it is weak lower semicontinuous, we note that $\mathscr{H}$ is convex by assumption. As the rest of the functional consists of a (non-)convex perturbation of lower order, Ioffe's theorem \cite[Theorem 7.5]{FonLeo2007MMCV} now gives the weak lower semicontinuity. Note that to apply Ioffe's theorem, we use the compact embedding $W^{2,p}(\Omega)\hookrightarrow W^{1,\infty}(\Omega)$, which follows from $p>d$ (comp.\ \ref{Assu:Hyperstress}). Finally, since the minimization is feasible as by assumption  $\calE_0(\chi_{k-1},c_{k-1})<\infty$, we obtain the existence of a minimizer $\chi_k$. In particular, we also have that $\calE_0(\chi_{k},c_{k-1})<\infty$.
\medskip

\textit{Diffusion step. } Given $\chi_k\in W^{2,p}(\Omega;\R^d)$ from the mechanical step and $c_{k-1}\in L^\infty(\Omega)$ from the previous time step, 
we now solve the diffusion equation \eqref{Eqn:TimeDisc:FiniteStrainChemPot:Reg} via a fixed-point argument, analogously to e.g.\ \cite{Jung2015BbEMfCDS, FHKM2022GEAE}. Let us fix $\wt\mu\in L^\infty(\Omega)$
and define $\widetilde c$ such that $\widetilde{\mu}\in \pl_c\Phi(\nabla\chi_k,\widetilde c)$ a.e.\ in $\Omega$ (see Remark~\ref{Rem:RelationMuC}). 
We solve the linear problem $a_{\widetilde{c}}(\mu,\psi)
= \langle\xi_{\widetilde{c}},\psi\rangle$ for all $\psi\in H^\theta(\Omega)$, 
where the bilinear and linear forms $a_{\widetilde{c}}:H^\theta(\Omega)\times H^\theta(\Omega)\to[0,\infty)$ and $\xi_{\widetilde{c}}\in H^\theta(\Omega)^*$
are given via
\begin{align*}
a_{\widetilde{c}}(\mu,\psi)&=\int_\Omega\!\Big\{  \calM(\nabla\chi_{k},\widetilde c)\nabla\mu \cdot \nabla\psi  + \eta\sum_{\abs{\beta}=\theta}\pl^\beta\mu\cdot \pl^\beta\psi \Big\}\dx + \int_{\pl\Omega}\!\!\kappa\mu\psi\dS,\\
\langle\xi_{\widetilde{c}},\psi\rangle &= -\int_{\Omega}\frac{\widetilde{c}-c_{k-1}}{\tau}\psi\dd x + \int_{\pl\Omega}\kappa\mu_{k,\ext}\psi\dd S.
\end{align*}
The existence of a unique solution $\mu\in H^\theta(\Omega)$
follows from the generalized Poincar\'e inequality for
the space $H^\theta(\Omega)$, see e.g.\ \cite[Ch.\ 2.1.4]{Tema1997IDDS}. Since $H^\theta(\Omega)$ embeds compactly
into $L^\infty(\Omega)$ (recall that $\theta> d/2)$, we can define
the map $\mathcal{S}: L^\infty(\Omega)\to L^\infty(\Omega)$, 
that maps $\widetilde\mu$ to $\mu$. We show that 
$\mathcal{S}$ is continuous, compact, and that for some $\lambda\in(0,1]$ the set 
$\mathsf{S}_\lambda:=\{\mu\in L^\infty(\Omega)\,|\, \mu = \lambda \mathcal{S}(\mu) \}$ is bounded.
Then, $\mathcal{S}$ has a fixed-point $\mu_k\in H^\theta(\Omega)\subset L^\infty(\Omega)$ by Schauder's fixed point theorem, see \cite[Thm.\ 11.3]{GilTru2001EPDE}.

To show the continuity of $\mathcal{S}$, consider 
a sequence of chemical potentials $\widetilde\mu^{(n)}\to \widetilde \mu$ in $L^\infty(\Omega)$ as $n\to \infty$ 
and the associated
densities $\widetilde{c}^{(n)}$. In particular,
we have that $\|\widetilde \mu^{(n)}\|_{L^\infty(\Omega)}\leq K$
for some $K>0$ independent of $n$ and thus also
$\|\widetilde{c}^{(n)}\|_{L^\infty(\Omega)} \leq C(K)$. Let now $\mu^{(n)} := \mathcal{S}(\widetilde\mu^{(n)})$.
We have to show that this sequence has a limit 
$\mu\in L^\infty(\Omega)$ such that $\mathcal{S}(\widetilde{\mu})=\mu$. Choosing $\psi = \mu^{(n)}$
in the equation for $\mu^{(n)}$ gives the uniform
estimate $\|\mu^{(n)}\|_{H^\theta(\Omega)} \leq C$ via standard arguments. Thus, we can find a (non-relabeled) subsequence 
and a limit $\mu\in H^\theta(\Omega)$
such that $\mu^{(n)}\wto \mu$ in $H^\theta(\Omega)$.
Moreover, by possibly passing to another subsequence
we can assume that $\widetilde{c}^{(n)}\to \widetilde{c}$
a.e.\ in $\Omega$ and by dominated convergence
also in $L^{\tilde p}(\Omega)$ for every $\tilde{p}\in [1,\infty)$. Thus, we can pass to the limit in 
$a_{\widetilde{c}^{(n)}}(\mu^{(n)},\psi) = \langle\xi_{\widetilde{c}^{(n)}},\psi\rangle$
to find that $\mu$ is the unique solution of
$a_{\widetilde{c}}(\mu,\psi) = \langle\xi_{\widetilde{c}},\psi\rangle$.
In particular, since it is unique, we get that $\mathcal{\widetilde{\mu}}= \mu$ and all converging subsequences converge to the same limit $\mu$. Thus, $\mathcal{S}$ is continuous.
The compactness of $\mathcal{S}$ follows similarly, as
the images of bounded sets in $L^\infty(\Omega)$  
under $\mathcal{S}$ are bounded in $H^\theta(\Omega)$.
Due to the compact embedding of $H^\theta(\Omega)$
into $L^\infty(\Omega)$ the claim follows.

Finally, we have to show that the set
$\mathsf{S}_\lambda$ for some $\lambda\in(0,1]$ is bounded. 
Let $\mu_\lambda\in\mathsf{S}_\lambda$, i.e., $\mu_\lambda$ satisfies
$a_{c_\lambda}(\mu_\lambda,\psi)= \lambda \langle\xi_{c_\lambda},\psi\rangle$
for all $\psi\in H^\theta(\Omega)$,
where $\mu_\lambda = \pl_c\Phi(\nabla\chi_k,c_\lambda)$.
Choosing $\psi = \mu_\lambda$ in the equation yields
\begin{multline*}
a_{c_\lambda}(\mu_\lambda,\mu_\lambda) = 
\int_\Omega\!\Big\{  \calM(\nabla\chi_{k}, c_\lambda)\nabla\mu_\lambda \cdot \nabla\mu_\lambda  + \eta\sum_{\abs{\beta}=\theta}|\pl^\beta\mu_\lambda|^2\Big\}\dx + \int_{\pl\Omega}\!\!\kappa\mu_\lambda^2\dS\\
= \langle\xi_{c_\lambda},\mu_\lambda\rangle=\lambda\int_{\pl\Omega}\kappa\mu_\lambda \mu_{k,\ext}\dd S -
\lambda\int_\Omega \frac{c_\lambda-c_{k-1}}{\tau}\mu_\lambda \dd x\\
\leq \lambda\int_{\pl\Omega}\kappa\mu_\lambda \mu_{k,\ext}\dd S 
+\frac{\lambda}{\tau}\int_\Omega \Phi(\nabla \chi_k,c_{k-1})\dd x - \frac{\lambda}{\tau}\int_{\Omega}\Phi(\nabla\chi_k,c_\lambda)\dd x.
\end{multline*}
Applying standard H\"older's and Young's inequalities for the boundary integral on the right-hand side, we conclude that
\begin{multline*}
\lambda \int_{\Omega}\Phi(\nabla \chi_k,c_k)+
\tau\int_\Omega\!\Big\{  \calM(\nabla\chi_{k}, c_\lambda)\nabla\mu_\lambda \cdot \nabla\mu_\lambda  + \eta\sum_{\abs{\beta}=\theta}|\pl^\beta\mu_\lambda|^2\Big\}\dx + \tau\kappa\big(1-\frac{\lambda}{2}\big)\int_{\pl\Omega}\!\!\mu_\lambda^2\dS\\
\leq \int_{\pl\Omega}\!\!\frac{\kappa\lambda}{2}\mu_{k,\ext}^2\dd S+ \lambda\int_\Omega \Phi(\nabla \chi_k,c_{k-1})\dd x,
\end{multline*}
and thus that $\|\mu_\lambda\|_{L^\infty(\Omega)}\leq C_{\infty,\theta}\|\mu_\lambda\|_{H^\theta(\Omega)} \leq C$
since $\Phi$ is bounded from below and the right-hand side is bounded.

Thus, Schauder's fixed-point theorem yields the existence of a fixed point $\mu_k\in H^\theta(\Omega)$. With $c_k$ defined as above, this fixed-point is a solution to the regularized diffusion equation
\eqref{Eqn:TimeDisc:FiniteStrainChemPot:Reg}. Moreover, we also
obtain from the estimate above for $\lambda=1$ 
that $\calE_0(\chi_k,c_k)<\infty$. Note that by Remark~\ref{Rem:RelationMuC}, we have that
$c_k\in L^\infty(\Omega)$ and $c_k>0$ a.e.\ in $\Omega$.

\end{proof}

To pass to the limit $(\tau,\eta) \to 0$, 
we now introduce the piecewise constant and piecewise affine interpolants
with respect to the time-discrete solutions obtained in Lemma~\ref{lemma:DiscreteSolutionsExistence}
by 
\begin{alignat}{3}
    &\chibar_{\eta,\tau}(0) = \chi_0^\eta, &&\chibar_{\eta,\tau}(t) = \chi_k^\eta && \text{for}\ t\in(t_{k-1},t_k], \\
    &\ubar\chi _{\eta,\tau}(0) = \chi_0^\eta, &&\ubar\chi_{\eta,\tau}(t) = \chi_{k-1}^\eta && \text{for}\ t\in[t_{k-1},t_k), \\
    &\chihat_{\eta,\tau}(0) = \chi_0^\eta,\hspace{10pt} 
        &&\chihat_{\eta,\tau}(t) = \frac{t-t_{k-1}}{\tau}\chi_k^\eta + \frac{t_k-t}{\tau}\chi_{k-1}^\eta\hspace{10pt} 
        && \text{for}\ t\in(t_{k-1},t_k], 
\end{alignat} 
	where, as before, $t_k=k\tau$ for $k=0,\ldots, N$. Similarly, we define the interpolants with respect to $c_k$, $\mu_k$, $\ell_k$,
	and $\mu_{k,\ext}$. 

	\begin{Rem}
		\label{Rem:TimeDisc:InterpolantEstimates}
		Note that the strong convergence of the piecewise affine interpolants together with a suitable estimate for the time derivative of the latter also implies the strong convergence of the piecewise constant interpolants, as for any Banach space $X$ the following estimates hold:%
		\begin{equation}
		\label{Eqn:TimeDisc:InterpolantEstimateL2}
		\norm{\wh \chi_\tau - \ubar\chi_\tau}_{L^2(0,T;X)} \leq \frac{\tau}{\sqrt{3}}\norm{\dot\chihat_\tau}_{L^2(0,T;X)},
		\end{equation}
		\begin{equation}
		\label{Eqn:TimeDisc:InterpolantEstimateLinfty}
		\norm{\wh \chi_\tau - \ubar\chi_\tau}_{L^\infty(0,T;X)} \leq \tau^{1/2}\norm{\dot\chihat_\tau}_{L^2(0,T;X)}.
		\end{equation}
		Of course, similar estimates also hold for $\chibar_\tau$.
	\end{Rem}
	
Using the piecewise affine and constant interpolants defined above, 
we can rewrite equation \eqref{Eqn:TimeDisc:FiniteStrainStress} for the deformation in the form
\begin{subequations}
\label{Eqn:TimeDisc:WeakFiniteStrainSystem}
    \begin{multline} 
    \label{Eqn:TimeDisc:WeakFiniteStrainStress}
    \int_0^T\int_\Omega \big(\sigma_\el(\nabla\chibar_{\eta,\tau},\ubar c_{\eta,\tau}) 
    + \sigma_\vi(\nabla\ubar\chi_{\eta,\tau},\nabla\dot\chihat_{\eta,\tau},\ubar c_{\eta,\tau})\big) : \nabla\phi\dd x\dd t\\  
    + \int_0^T\int_\Omega\mfh(\rmD^2 \chibar_{\eta,\tau}) \tdots \rmD^2\phi \dx\dt = \int_0^T\ip*{\wb \ell_\tau,\phi}\dt,
    \end{multline}
    for $\phi\in L^2(0,T;W_0^{2,p}(\Omega;\R^d))$ %
    and \eqref{Eqn:TimeDisc:FiniteStrainChemPot:Reg} for the concentration as
    \begin{multline}				\label{Eqn:TimeDisc:WeakFiniteStrainChemPot}
    \int_0^T \ip{\dot\chat_{\eta,\tau},\psi}\dt 
    + \int_0^T\int_\Omega \calM(\nabla\chibar_{\eta,\tau},\cbar_{\eta,\tau})\nabla\mubar_{\eta,\tau} \cdot \nabla\psi 
    + \eta \sum_{\abs{\beta}=\theta}\pl^\beta\mubar_{\eta,\tau} \cdot \pl^\beta\psi\dx \dt\\
    + \int_0^T\int_{\pl\Omega}\kappa(\mubar_{\eta,\tau} - \mubar_{\tau,\ext})\psi\dS\dt = 0,
    \end{multline}
    for $\psi\in L^{s'}(0,T;H^\theta(\Omega))$, where $\mubar_{\eta,\tau} \in \partial_{c}\Phi(\chibar_{\eta,\tau},\cbar_{\eta,\tau})$  a.e.\ in $\Omega$.
\end{subequations}	
		
We start by showing an energy-dissipation inequality.

\begin{Lemma}[Energy-dissipation inequality]
    Let $\chibar_{\eta,\tau}$, $\cbar_{\eta,\tau}$, $\ubar\chi_{\eta,\tau}$, $\ubar \cbar_{\eta,\tau}$, $\mubar_{\eta,\tau}$, and $\chihat_{\eta,\tau}$ denote the piecewise constant and affine
interpolants with respect to the time-discrete solutions
obtained in Lemma \ref{lemma:DiscreteSolutionsExistence}.
For $\calE_\tau(t,\chi,c) := \int_\Omega \Phi(\nabla\chi,c)+\mathscr{H}(\rmD^2\chi)\dx
        -\ip{\bar\ell_\tau(t),\chi}$,
the following time-discrete energy-dissipation inequality holds for $t_k = k\tau$, $k=1,\ldots,N$:
\begin{multline}\label{Eqn:TimeDiscEDIneq}
    \calE_\tau(t_k,\chibar_{\eta,\tau}(t_k),\cbar_{\eta,\tau}(t_k)) + 
	\int_0^{t_k} \int_\Omega\Big\{ \calM(\nabla\chibar_{\eta,\tau},\cbar_{\eta,\tau})\nabla\mubar_{\eta,\tau}{\cdot}\nabla\mubar_{\eta,\tau}
    +\eta\sum_{\abs{\beta}=\theta}\abs{\pl^\beta\mubar_{\eta,\tau}}^2\Big\}\dx\dt \\
			+ \int_0^{t_k}\int_{\pl\Omega}\kappa\mubar_{\eta,\tau}^2\dS\dt
+\int_0^{t_k}\!\calR\big(\ubar\chi_{\eta,\tau},\dot\chihat_{\eta,\tau},\ubar \cbar_{\eta,\tau}\big)\dt \\
\leq \calE_\tau(0,\chi_{0}^\eta,c_{0}^\eta) + \int_0^{t_k}\int_{\pl\Omega}\kappa\mubar_{\eta,\tau}\mubar_{\tau,\ext}\dS\dt 
- \int_0^{t_k}\ip{\dot {\widehat{\ell}}_\tau,\ubar\chi_{\eta,\tau}}\dt.
\end{multline}
\end{Lemma}
\begin{proof}
    We drop the index $\eta$ throughout the proof.
    
Using the definition of the subdifferential, we get with $\mu_k = \pl_c\Phi(\nabla\chi_k,c_k)$ that a.e.\ in $\Omega$
\[
	\Phi(\nabla\chi_k,c_{k-1})\geq \Phi(\nabla\chi_k,c_k)+\mu_k(c_{k-1}-c_k). %
\] 
Thus, it follows for all $k=1,\ldots,N$, 
after integration over $\Omega$ 
(and adding the hyperstress potential 
$\mathscr{H}$ to both sides) that
\begin{equation*}
\int_\Omega \Phi(\nabla\chi_k,c_k)
+\mathscr{H}(\rmD^2\chi_k)\dx 
- \int_\Omega \Phi(\nabla\chi_k,c_{k-1})+\mathscr{H}(\rmD^2\chi_k)\dx -\int_\Omega \mu_k(c_k {-} c_{k-1})\dx 	\leq 0.
\end{equation*} 
Using the time-discrete diffusion equation in
\eqref{Eqn:TimeDisc:FiniteStrainChemPot:Reg} 
with test function $\psi = \mu_k\in H^\theta(\Omega)$ for the
last term on the left-hand side, we now obtain the discrete
energy inequality
\begin{multline}
    \label{Eqn:TimeDisc:CalcAPriori1}
			\mathcal{E}_0(\chi_k,c_k) - \mathcal{E}_0(\chi_k,c_{k-1}) +\tau\int_\Omega \calM(\nabla\chi_k,c_k)\nabla\mu_k\cdot\nabla\mu_k\dx 
			\\+ \tau\int_{\pl\Omega}\kappa(\mu_k^2 - \mu_k\mu_{k,\ext})\dS + \tau\eta\int_\Omega\sum_{\abs{\beta}=\theta}\abs{\pl^\beta\mu_k}^2\dd x %
			\leq 0,
\end{multline}
where, as before, $\calE_0(\chi,c)=\int_\Omega \Phi(\nabla\chi,c)+\mathscr{H}(\rmD^2 \chi)\dd x$.

Next, we establish the analogous estimate for the mechanical step. This estimate, however, follows directly from  $\chi_k$ being the solution of the minimization problem \eqref{Eqn:TimeDisc:MinimizerChi}. Indeed, choosing $\chi_{k-1}$ as competitor in the minimization gives 
\begin{equation}\label{Eqn:EstimateMinimizationChi}
\calE_0(\chi_k,c_{k-1})  +
	\tau\calR(\chi_{k-1},\delta_\tau\chi_k,c_{k-1}) -\ip{\ell_k,\chi_k}
	\leq \calE_0(\chi_{k-1},c_{k-1})-\ip{\ell_k,\chi_{k-1}}.
\end{equation}
To obtain the estimate in \eqref{Eqn:TimeDiscEDIneq}, we use the estimate \eqref{Eqn:EstimateMinimizationChi} for $\calE_0(\chi_k,c_{k-1})$ in the estimate \eqref{Eqn:TimeDisc:CalcAPriori1} to get
\begin{multline*}
\mathcal{E}_0(\chi_i,c_i)  +\tau\int_\Omega\Big\{ \calM(\nabla\chi_{i},c_{i})\nabla\mu_i\cdot\nabla\mu_i
+\eta\sum_{\abs{\beta}=\theta}\abs{\pl^\beta\mu_i}^2
\Big\}\dd x
\\
+ \tau\int_{\pl\Omega}\kappa(\mu_i^2 - \mu_i\mu_{i,\ext})\dS
+	\tau\calR(\chi_{i-1},\delta_\tau\chi_i,c_{i-1}) 
	\leq \calE_0(\chi_{i-1},c_{i-1})+\ip{\ell_i,\chi_i{-}\chi_{i-1}}.
\end{multline*}
The inequality in \eqref{Eqn:TimeDiscEDIneq} then follows from summing over $i=1,\ldots, k$, summation by parts for the loading part, and the definition of the piecewise constant and affine interpolants. In particular, the summation-by-parts formula reads
 \begin{align*}
\sum_{i=1}^k \tau\ip{\ell_i,\delta_\tau \chi_i} &= \ip{\ell_k,\chi_k}-\ip{\ell_0,\chi_0} - \sum_{i=1}^{k}\tau \ip*{\delta_\tau \ell_{i},\chi_{i-1}}.
\end{align*}
\end{proof}    

The time-discrete energy-dissipation inequality gives rise to uniform a priori bounds (with respect to $\tau$ and $\eta$)
for the interpolants of the time-discrete solutions $(\chi_k,\mu_k,c_k)$ obtained in Lemma~\ref{lemma:DiscreteSolutionsExistence}.
We assume in the following that the initial conditions satisfy 
\begin{equation}\label{Eqn:ConvInitialCond}
\chi_0^\eta \to \chi_0\text{ in }W_\id^{2,p}(\Omega;\R^d),\quad c_0^\eta \to c_0\text{ in }L^1(\Omega),\quad\text{and}\quad\calE_0(\chi_0^\eta,c_0^\eta) \to \calE_0(\chi_0,c_0)
\end{equation}
 as $\eta\to 0$.
 
\begin{Lemma}[A priori estimates]
\label{lem:AprioriEstimatesTimeDisc}
Assuming that \ref{Assu:Hyperstress}--\ref{Assu:RegularityInit}, as well as \eqref{Eqn:ConvInitialCond} are satisfied, there exists $C>0$ (independent of $\tau$ and $\eta$) such that for $\tau, \eta>0$ small enough, the following uniform estimates are satisfied
\begin{subequations}
\begin{align}
\label{Eqn:TimeDisc:APrioriChi}
&\norm{\chibar_{\eta,\tau}}_{L^\infty(0,T;W^{2,p}_\id(\Omega))} 
    + \norm{\dot{\chihat}_{\eta,\tau}}_{L^2(0,T;H^1_0(\Omega))} 
    + \|(\det\nabla\chibar_{\eta,\tau})^{-1}\|_{L^\infty(0,T;L^q(\Omega))} \leq C,\\
\label{Eqn:TimeDisc:APrioriMu}
    &\sqrt{\eta}\norm{\mubar_{\eta,\tau}}_{L^2(0,T;H^\theta(\Omega))} 
        + \norm{\sqrt{\kappa}\mubar_{\eta,\tau}}_{L^2([0,T]\times\pl\Omega)} \leq C,\\
\label{Eqn:TimeDisc:APrioriC}
    &\norm{\cbar_{\eta,\tau}}_{L^\infty(0,T;L\log L(\Omega))} 
        + \big\|\nabla\cbar^{\frac{m}{2}}_{\eta,\tau}\big\|_{L^2(0,T;L^2(\Omega))} 
        + \|\dot{\chat}_{\eta,\tau}\|_{L^s(0,T;H^\theta(\Omega)^*)} \leq C,\\
\label{Eqn:TimeDisc:APrioriMob}
    & \norm{\calM(\nabla\chibar_{\eta,\tau},\cbar_{\eta,\tau})\nabla\mubar_{\eta,\tau}}_{L^s([0,T]\times \Omega)} 
\leq C_\eta.
\end{align}
where $s=\frac{md+2}{md+1}>1$.

If additionally $\gamma_1 >0$ in Assumption \ref{Assu:FreeEn}(ii) then 
we also have that
\begin{align}
   \label{Eqn:APriori:Reg:c:CaseII}
    &\norm{\cbar_{\eta,\tau}}_{L^\infty(0,T;L^{2+r}(\Omega))} 
    + \sum_{\omega\in \{0,\frac{1+r}{2},\,1+r\}} \big\|{\nabla \cbar^{\frac{m}{2}+\omega}_{\eta,\tau}}\big\|^2_{L^2(0,T;L^2(\Omega))}
    + \norm{\dot \chat_{\eta,\tau}}_{L^s(0,T;H^\theta(\Omega)^*)} 
    \leq C,\\
   \label{Eqn:APriori:Reg:mob:CaseII}
    &\norm*{\calM(\nabla\chi_\eta, c_\eta)\nabla
    \mu_\eta}_{L^s([0,T]\times \Omega)}\leq C,
\end{align}
where $s=\min\{\frac{md+2(r+2)}{md+r+2},\frac{d(m+r+1)+2(r+2)}{d(m+r+1)+r+2}\}>1$.%
\end{subequations}
\end{Lemma}

\begin{proof} 
\emph{Step 1.} To obtain the a priori estimate for
$\chibar_{\eta,\tau}$, we exploit the energy-dissipation estimate in \eqref{Eqn:TimeDiscEDIneq}. We absorb the boundary integral on the right-hand side via Young's inequality in the corresponding boundary terms on the left-hand side. Moreover, neglecting any nonnegative terms on the left-hand side gives
\begin{multline}
    \calE_\tau(t_k,\chibar_{\eta,\tau}(t_k), \cbar_{\eta,\tau}(t_k))
    \leq \calE_\tau(0,\chi_0,c_0) + \frac{1}{2}\int_0^{t_k} 
\|\sqrt{\kappa}\mubar_{\tau,\ext}\|^2_{L^2(\partial\Omega)}\dt \\
+ \int_0^{t_k} \big\|\dot {\widehat{\ell}}_\tau\big\|_{W^{2,p}(\Omega;\R^d)^*}
    \norm{\ubar\chi_{\eta,\tau}}_{W^{2,p}_\id(\Omega;\R^d)}\dt
    \leq C_0 + \Lambda_1\int_0^{t_k}\|\ubar\chi_{\eta,\tau}(t)\|_{W^{2,p}_\id(\Omega;\R^d)}\dt,
\end{multline}
where we used Assumptions~\ref{Assu:Permeability} and \ref{Assu:Loading} for the data $\mu_\ext$ and $\ell$, respectively, and set $\Lambda_1:=\|\dot\ell\|_{L^\infty(0,T;W^{2,p}(\Omega;\R^d)^*)}$. However, due to Assumptions \ref{Assu:Hyperstress}, \ref{Assu:FreeEn}(i), and \ref{Assu:Loading}, we find a constant $C_\calE>0$ such that the energy $\calE_\tau$ satisfies the lower estimate
\begin{equation}\label{Eqn:LowerEstimateEnergy}
\calE_\tau(t,\chi,c)\geq C_\calE\|\chi\|_{W^{2,p}_\id(\Omega;\R^d)} + C_{\Phi,0}\|(\det\nabla\chi)^{-1}\|_{L^q(\Omega)}^q-1/C_\calE.
\end{equation}
Thus, setting $e_k:=\calE_\tau(t_k,\chibar_{\eta,\tau}(t_k),\cbar_{\eta,\tau}(t_k))$, we have the time-discrete estimate 
\[
e_k \leq C_0+ \frac{\Lambda_1 T}{C^2_\calE} + \frac{\Lambda_1}{C_\calE}\sum_{i=0}^{k-1}  \tau e_i,\quad k =1,\ldots,N.
\]

An application of the discrete Gronwall lemma (see e.g.~\cite[Lemma 1.4.2]{QuaVal1994NAPD}) now gives the bound 
\begin{equation}\label{Eqn:GrownwallTimeDiscrete}
e_k \leq \big(C_0+\frac{\Lambda_1 T}{C_\calE^2}\big)\exp\Big(\frac{\Lambda_1}{C_\calE}k\tau\Big).
\end{equation}
Hence, estimate \eqref{Eqn:LowerEstimateEnergy} gives the bound for 
$\chibar_{\eta,\tau}$ in $L^\infty(0,T;W^{2,p}_\id(\Omega;\R^d))$ (cf.~\eqref{Eqn:TimeDisc:APrioriChi}).

Note that as an immediate consequence of \eqref{Eqn:LowerEstimateEnergy}  and  \eqref{Eqn:GrownwallTimeDiscrete}, 
we also obtain that $1/(\det\nabla\chibar_{\eta,\tau})$ is uniformly
bounded in $L^\infty(0,T;L^q(\Omega))$ for $q\geq pd/(p{-}d)$
(comp.\ Assumption~\ref{Assu:FreeEn}(i)). Consequently, the
Healey--Kr\"omer lemma \cite[Thm.~3.1]{MieRou2020TVKVR} can be
applied, which gives a uniform constant $C_{\text{HK}}>0$ and the
lower bound 
\[
\det\nabla\chibar_{\eta,\tau}(t,x) \geq C_\text{HK}>0\quad\text{for all }(t,x)\in[0,T]\times\Omega.
\] 
(Note that by assumption, also $\det\nabla \chi_0(x)\geq \rho_0>0$.)
In particular, we see that for some $R>0$ large enough, 
$\nabla\chibar_{\eta,\tau}(t,x) \in \mathsf{F}_R$ for all $(t,x)\in[0,T]\times\overline\Omega$. \medskip

\emph{Step 2.} To obtain the a priori estimate for $\dot\chihat_{\eta,\tau}$ in $L^2(0,T;H^1(\Omega))$, we use assumption \ref{Assu:ViscousStress} to estimate
\begin{align*}
\int_0^T\calR(\ubar\chi_{\eta,\tau},\dot\chihat_{\eta,\tau},&\ubar c_{\eta,\tau})\dt
    = \int_0^T\int_\Omega \zeta(\nabla\ubar\chi_{\eta,\tau},\nabla\dot\chihat_{\eta,\tau}, \ubar c_{\eta,\tau})\dx\dt\\
    &= \int_0^T\int_\Omega \wh\zeta\big((\nabla\ubar\chi_{\eta,\tau})^\top \nabla\ubar\chi_{\eta,\tau}, (\nabla\ubar\chi_{\eta,\tau})^\top \nabla\dot\chihat_{\eta,\tau} 
    + (\nabla\dot\chihat_{\eta,\tau})^\top \nabla\ubar\chi_{\eta,\tau}, \ubar c_{\eta,\tau}\big)\dx\dt\\
    &\geq C_{\zeta,1}\int_0^T \int_\Omega \abs*{(\nabla\ubar\chi_{\eta,\tau})^\top \nabla\dot\chihat_{\eta,\tau} 
+ (\nabla\dot\chihat_{\eta,\tau})^\top \nabla\ubar\chi_{\eta,\tau}}^2\dd x\dt\\
    &\geq C_{\zeta,1}\int_0^T \norm{\dot\chihat_{\eta,\tau}}^2_{H^1(\Omega)}\dt 
= C_{\zeta,1}\norm{\dot\chihat_{\eta,\tau}}^2_{L^2(0,T;H^1(\Omega))},
\end{align*}
		where the last inequality follows from the generalized Korn's inequality as in \cite[Cor.~3.4]{MieRou2020TVKVR}. %
        The uniform boundedness of the left-hand side now follows directly from the energy-dissipation inequality \eqref{Eqn:TimeDiscEDIneq} and hence the proof of \eqref{Eqn:TimeDisc:APrioriChi} is complete.
		\medskip

		\emph{Step 3.} To obtain the a priori estimates for $\mubar_{\eta,\tau}$, 
we again exploit the energy-dissipation inequality~\eqref{Eqn:TimeDiscEDIneq}  
as well as the generalized Poincar\'{e} inequality (see \cite[Section\ 2.1.4]{Tema1997IDDS}) to obtain the estimate
\[
\sqrt{\eta}\norm{\mubar_{\eta,\tau}}_{L^2(0,T;H^\theta(\Omega))} + \norm{\sqrt{\kappa}\mubar_{\eta,\tau}}_{L^2([0,T]\times\pl\Omega)}\leq C
\]
for some constant $C>0$ independent of $\tau$ and $\eta$. Thus, we have shown~\eqref{Eqn:TimeDisc:APrioriMu}.

		\medskip

        \emph{Step 4.} The first estimate of \eqref{Eqn:TimeDisc:APrioriC}, i.e., the estimate for $\cbar_{\eta,\tau}$ in $L^\infty(0,T;L\log L(\Omega))$ %
follows directly from the 
$L^\infty(0,T;L^\infty(\Omega;\R^{d\times d}))$-bound for $\nabla\chibar_{\eta,\tau}$, Lemma \ref{Lemma:BelowEstimatePhiMix}, and the energy-dissipation inequality \eqref{Eqn:TimeDiscEDIneq}. Similarly, when $\gamma_1>0$, we also obtain the bound for $\cbar_{\eta,\tau}$ in $L^\infty(0,T;L^{2+r}(\Omega))$.  

To obtain the bounds for the gradient terms, %
we now use that
    \begin{align*}
        \nabla\mubar_{\eta,\tau} &= \pl^2_{Fc}\Phi(\nabla\chibar_{\eta,\tau},\cbar_{\eta,\tau})\rmD^2\chibar_{\eta,\tau} + \pl^2_{cc}\Phi(\nabla\chibar_{\eta,\tau},\cbar_{\eta,\tau})\nabla \cbar_{\eta,\tau},
    \end{align*}
    which by the binomial formula implies that
    \begin{multline*}
        \int_0^T\int_\Omega \calM(\nabla\chibar_{\eta,\tau},\cbar_{\eta,\tau})\nabla\mubar_{\eta,\tau}\cdot\nabla\mubar_{\eta,\tau} \dx\dt 
\\
= \int_0^T\int_\Omega \calM(\nabla\chibar_{\eta,\tau},\cbar_{\eta,\tau})\Big\{\abs*{\pl^2_{Fc}\Phi(\nabla\chibar_{\eta,\tau},\cbar_{\eta,\tau})\rmD^2\chibar_{\eta,\tau}}^2%
        + \abs*{\pl^2_{cc}\Phi(\nabla\chibar_{\eta,\tau},\cbar_{\eta,\tau})\nabla \cbar_{\eta,\tau}}^2\\+ 2(\pl^2_{Fc}\Phi(\nabla\chibar_{\eta,\tau},\cbar_{\eta,\tau})\rmD^2\chibar_{\eta,\tau}) \cdot \pl^2_{cc}\Phi(\nabla\chibar_{\eta,\tau},\cbar_{\eta,\tau})\nabla \cbar_{\eta,\tau}\Big\}\dx\dt. %
    \end{multline*}
    Thus, using Young's inequality with $\epsilon>0$ for the last term in the brackets, 
    we see that 
\begin{multline}
\int_0^T\int_\Omega \calM(\nabla\chibar_{\eta,\tau},\cbar_{\eta,\tau})\nabla\mubar_{\eta,\tau}\cdot\nabla\mubar_{\eta,\tau} \dx\dt\\ \geq C\int_0^T\int_\Omega \pl^2_{cc}\Phi(\nabla\chibar_{\eta,\tau},\cbar_{\eta,\tau})^2\calM(\nabla\chibar_{\eta,\tau},\cbar_{\eta,\tau})\nabla \cbar_{\eta,\tau}\cdot\nabla \cbar_{\eta,\tau} \dx\dt \\
- C\int_0^T\int_\Omega \abs{\calM(\nabla\chibar_{\eta,\tau},\cbar_{\eta,\tau})}\abs*{\pl^2_{Fc}\Phi(\nabla\chibar_{\eta,\tau},\cbar_{\eta,\tau})\rmD^2\chibar_{\eta,\tau}}^2 \dx\dt. \label{Eqn:FiniteStrain:Interm1}
    \end{multline}\smallskip

    We now discuss the two cases from Assumption~\ref{Assu:FreeEn}(ii), separately, 
    i.e., $\gamma_1=\gamma_2=0$ or $0<\gamma_1\leq \gamma_2$.
    
    \underline{\textbf{Case I}: $\gamma_1=\gamma_2=0$.} 
To estimate the first integral on the right-hand side of \eqref{Eqn:FiniteStrain:Interm1}, we use that $\nabla\chibar_{\eta,\tau} \in \mathsf{F}_R$ and Assumptions \ref{Assu:MobilityTensor}, \ref{Assu:FreeEn}(ii) to obtain
    \begin{multline}
    \int_0^T\int_\Omega  \pl^2_{cc}\Phi(\nabla\chibar_{\eta,\tau},\cbar_{\eta,\tau})^2\calM(\nabla\chibar_{\eta,\tau},\cbar_{\eta,\tau})\nabla \cbar_{\eta,\tau}\cdot\nabla \cbar_{\eta,\tau} \dx\dt \\\geq C\int_0^T\int_\Omega \cbar_{\eta,\tau}^{m-2} \abs{\nabla \cbar_{\eta,\tau}}^2\dx\dt
    \geq C\big\|{\nabla \cbar_{\eta,\tau}^{\frac{m}{2}}}\big\|_{L^2(0,T;L^2(\Omega))} .
    \label{Eqn:FiniteStrain:Interm2}
    \end{multline}

For the second integral on the right-hand side of \eqref{Eqn:FiniteStrain:Interm1}, we note that by \ref{Assu:FreeEn}(iii), Hölder's inequality, and the uniform bound for $\rmD^2\chibar_{\eta,\tau}$ in $L^\infty(0,T;L^p(\Omega;\R^{d\times d\times d}))$
\begin{multline*}
\int_0^T\int_\Omega \abs{\calM(\nabla\chibar_{\eta,\tau},\cbar_{\eta,\tau})}
\abs*{\pl^2_{Fc}\Phi(\nabla\chibar_{\eta,\tau},\cbar_{\eta,\tau})\rmD^2\chibar_{\eta,\tau}}^2 \dx\dt 
\\\leq \int_0^T\int_\Omega \cbar_{\eta,\tau}^{m+2\alpha}\abs{\rmD^2\chibar_{\eta,\tau}}^2\dx\dt
\leq C\norm{\cbar^{m+2\alpha}_{\eta,\tau}}_{L^1(0,T;L^\frac{p}{p-2}(\Omega))}.
\end{multline*}
We now distinguish two cases: First, if $0\leq \frac{p}{p-2}(m+2\alpha)\leq 1$, a uniform upper estimate for the right-hand side follows immediately from the $L^\infty(0,T;L^1(\Omega))$-bound for $\cbar_{\eta,\tau}$ (cf.~\eqref{Eqn:TimeDisc:APrioriC}). 
For the second case, when $\frac{p}{p-2}(m+2\alpha)>1$, we note that $\frac{p}{p-2}\leq \frac{m}{2}\frac{2d}{d-2}\frac{1}{m+2\alpha} =: \omega_2\in (1,\infty)$, which follows from Assumption~\ref{Assu:FreeEn}(iii) and $p>d$. Since $\alpha <0$ in this case,
we also have that $1<\frac{m}{2}\frac{2}{m+2\alpha}=:\omega_1$. 
Thus, we now find that
    \begin{align*}
        \norm*{\cbar^{m+2\alpha}_{\eta,\tau}}_{L^1(0,T;L^\frac{p}{p-2}(\Omega))} 
	    \leq C\norm*{\cbar_{\eta,\tau}^{m+2\alpha}}_{L^{\omega_1}(0,T;L^{\omega_2}(\Omega))}
        = C\big\|{\cbar_{\eta,\tau}^{\frac{m}{2}}}\big\|^{2/\omega_1}_{L^2(0,T;L^{\frac{2d}{d-2}}(\Omega))}.
	\end{align*}
    Next, we note that $\frac{2}{\omega_1}< 2$, and thus we can use Young's inequality with $\epsilon>0$ to obtain
    \begin{align*}
    \norm*{\cbar^{m+2\alpha}_{\eta,\tau}}_{L^1(0,T;L^\frac{p}{p-2}(\Omega))}
        &\leq  C(\epsilon) + \epsilon C\big\|{\cbar_{\eta,\tau}^{\frac{m}{2}}}\big\|^2_{L^2(0,T;L^{\frac{2d}{d-2}}(\Omega))}.
	\end{align*}
    Using the Sobolev embedding $H^1(\Omega)\hookrightarrow L^{\frac{2d}{d-2}}(\Omega)$, we obtain
    \begin{align*}
    	\int_0^T\int_\Omega \abs{\calM(\nabla\chibar_{\eta,\tau},\cbar_{\eta,\tau})}\abs*{\pl^2_{Fc}\Phi(\nabla\chibar_{\eta,\tau},\cbar_{\eta,\tau})\rmD^2\chibar_{\eta,\tau}}^2\dx\dt &\leq C(\epsilon) + \epsilon C\big\|{\cbar_{\eta,\tau}^{\frac{m}{2}}}\big\|^2_{L^2(0,T;H^1(\Omega))}.
    \end{align*}
    Finally, we note that $\cbar_{\eta,\tau}^{\frac{m}{2}} $ is uniformly bounded in 
$L^2(0,T;L^1(\Omega))$ due to the boundedness of the energy and 
$m \leq 2$. 
Thus we can use a Poincar\'e-type inequality to find
    \begin{equation}
    \label{Eqn:FiniteStrain:Interm3}
    	\int_0^T\int_\Omega \abs{\calM(\nabla\chibar_{\eta,\tau},\cbar_{\eta,\tau})}\abs*{\pl^2_{Fc}\Phi(\nabla\chibar_{\eta,\tau},\cbar_{\eta,\tau})\rmD^2\chibar_{\eta,\tau}}^2\dx\dt \leq C(\epsilon) + \epsilon C\big\|{\nabla \cbar_{\eta,\tau}^{\frac{m}{2}}}\big\|^2_{L^2(0,T;L^2(\Omega))}.
    \end{equation}
	Choosing $\epsilon>0$ sufficiently small, we can combine \eqref{Eqn:FiniteStrain:Interm1}, \eqref{Eqn:FiniteStrain:Interm2}, and \eqref{Eqn:FiniteStrain:Interm3} to obtain    
	\begin{equation*}
	\int_0^T\int_\Omega \calM(\nabla\chibar_{\eta,\tau},\cbar_{\eta,\tau})\nabla\mubar_{\eta,\tau}\cdot\nabla\mubar_{\eta,\tau} \dx\dt \geq C\big(\big\|{\nabla \cbar_{\eta,\tau}^{\frac{m}{2}}}\big\|_{L^2(0,T;L^2(\Omega))}^2 - 1\big).
    \end{equation*}    
    Thus, the energy-dissipation inequality \eqref{Eqn:TimeDiscEDIneq} now gives the uniform bound for $\nabla \cbar_{\eta,\tau}^\frac{m}{2}$ in $L^2(0,T;L^2(\Omega))$.
    \medskip

\underline{\textbf{Case II}: $\gamma_2\geq \gamma_1>0$.} To estimate the first integral, we again use that $\nabla\chibar_{\eta,\tau} \in \mathsf{F}_R$ and Assumptions \ref{Assu:MobilityTensor}, \ref{Assu:FreeEn}(ii) to obtain
\begin{multline}
    \int_0^T\int_\Omega \pl^2_{cc}\Phi(\nabla\chibar_{\eta,\tau},\cbar_{\eta,\tau})^2\calM(\nabla\chibar_{\eta,\tau},\cbar_{\eta,\tau}) \nabla \cbar_{\eta,\tau}\cdot \nabla \cbar_{\eta,\tau} \dx\dt 
    \\\geq C\int_0^T\int_\Omega \cbar_{\eta,\tau}^{m} (\cbar_{\eta,\tau}^{-1} +  \cbar_{\eta,\tau}^r)^2 \abs{\nabla \cbar_{\eta,\tau}}^2\dx\dt\geq C\sum_{i\in \{0,\frac{1+r}{2},\, 1+r\}} \big\|{\nabla c^{\frac{m}{2}+i}_{\eta,\tau}}\big\|^2_{L^2(0,T;L^2(\Omega))}.
    \label{Eqn:FiniteStrain:Interm2:CaseII}
\end{multline}
For the second integral, we proceed as in Case I to obtain
\begin{align*}
\int_0^T\int_\Omega \abs{\calM(\nabla\chibar_{\eta,\tau},\cbar_{\eta,\tau})}
\abs*{\pl^2_{Fc}\Phi(\nabla\chibar_{\eta,\tau},\cbar_{\eta,\tau})\rmD^2\chibar_{\eta,\tau}}^2 \dx\dt 
\leq C\norm{\cbar^{m+2\alpha}_{\eta,\tau}}_{L^1(0,T;L^\frac{p}{p-2}(\Omega))}.
\end{align*}
Again, we consider two cases: First, if $0\leq \frac{p}{p-2}(m+2\alpha)\leq 2+r$ it follows from the $L^\infty(0,T;L^{2+r}(\Omega))$-bound for $\cbar_{\eta,\tau}$ that the left-hand side is uniformly bounded. If $\frac{p}{p-2}(m+2\alpha)>2+r$, we define $\omega_0$ by $\omega_0 = \frac{m+1+r}{2}$ in Case IIa, and $\omega_0 = \frac{m}{2}{+}1{+}r$ in Case IIb. Next, we note that $\frac{p}{p-2}\leq \frac{2d}{d-2}\frac{\omega_0}{m+2\alpha} =: \omega_2$, which follows from Assumption~\ref{Assu:FreeEn}(iii) and $p>d$. Also, $1< \frac{2\omega_0}{m+2\alpha} =: \omega_1$ so that
 \begin{align*}
        \norm*{\cbar^{m+2\alpha}_{\eta,\tau}}_{L^1(0,T;L^\frac{p}{p-2}(\Omega))} 
	    \leq C\norm*{\cbar_{\eta,\tau}^{m+2\alpha}}_{L^{\omega_1}(0,T;L^{\omega_2}(\Omega))}
        \leq C\big\|{\cbar_{\eta,\tau}^{\omega_0}}\big\|_{L^2(0,T;L^{\frac{2d}{d-2}}(\Omega))}^{2/\omega_1}
	\end{align*}
Since $\frac{2}{\omega_1}< 2$ we can proceed as in Case I and use Young's inequality with $\eps$, the Sobolev embedding $H^1(\Omega)\hookrightarrow L^{\frac{2d}{d-2}}(\Omega)$ and Poincar\'e's inequality (which is applicable since $\cbar_{\eta,\tau}^{\omega_0}\in L^2(0,T;L^1(\Omega))$ as $m\leq 3+r$ (Case IIa) or $m\leq 2$ (Case IIb)) to obtain
\begin{equation}
    \label{Eqn:FiniteStrain:Interm3:CaseII}
    	\int_0^T\int_\Omega \abs{\calM(\nabla\chibar_{\eta,\tau},\cbar_{\eta,\tau})}\abs*{\pl^2_{Fc}\Phi(\nabla\chibar_{\eta,\tau},\cbar_{\eta,\tau})\rmD^2\chibar_{\eta,\tau}}^2\dx\dt \leq C(\epsilon) + \epsilon C\norm*{\nabla \cbar_{\eta,\tau}^{\omega_0}}^2_{L^2(0,T;L^2(\Omega))}.
    \end{equation}
Thus, choosing $\eps$ sufficiently small and combining \eqref{Eqn:FiniteStrain:Interm1}, \eqref{Eqn:FiniteStrain:Interm2:CaseII} and \eqref{Eqn:FiniteStrain:Interm3:CaseII}, we obtain
\begin{equation*}
	\int_0^T\int_\Omega \calM(\nabla\chibar_{\eta,\tau},\cbar_{\eta,\tau})\nabla\mubar_{\eta,\tau}\cdot\nabla\mubar_{\eta,\tau} \dx\dt \geq C\Big(\sum_{i\in \{0,\frac{1+r}{2},\, 1+r\}} \big\|{\nabla c^{\frac{m}{2}+i}_{\eta,\tau}}\big\|^2_{L^2(0,T;L^2(\Omega))} - 1\Big).
    \end{equation*} 
The energy-dissipation inequality \eqref{Eqn:TimeDiscEDIneq} now gives the uniform bound for $\nabla \cbar_{\eta,\tau}^\frac{m}{2}$, $\nabla \cbar_{\eta,\tau}^\frac{m+1+r}{2}$ and $\nabla \cbar_{\eta,\tau}^{\frac{m}{2}+1+r}$ in $L^2(0,T;L^2(\Omega))$. \medskip

 \emph{Step 5.} 
 We show that the flux $\calM(\nabla\chibar_{\eta,\tau},\cbar_{\eta,\tau})\nabla\mubar_{\eta,\tau}$ is uniformly
 bounded in $L^s([0,T]\times\Omega)$ (cf.~\eqref{Eqn:TimeDisc:APrioriMob}). 
We distinguish again between the Cases I and II in Assumption~\ref{Assu:FreeEn}(iii).
 
 \underline{\textbf{Case I}: $\gamma_1=\gamma_2=0$.} We estimate for $1<s=\frac{md+2}{md+1}<2$
\begin{align*}
    \norm{\calM(\nabla\chibar_{\eta,\tau},\cbar_{\eta,\tau})\nabla\mubar_{\eta,\tau}}_{L^s([0,T]\times\Omega)}^s 
	&\leq C\int_0^T\int_\Omega \big(\cbar^m_{\eta,\tau}\big|\pl_{Fc}^2\Phi \rmD^2\chibar_{\eta,\tau} 
	+ \pl_{cc}^2\Phi\nabla \cbar_{\eta,\tau}\big|\big)^s\dx\dt\\
    &\leq C\int_0^T\int_\Omega \big(\cbar_{\eta,\tau}^{s(m+\alpha)}\abs{\rmD^2\chibar_{\eta,\tau}}^s + \cbar_{\eta,\tau}^{s(m-1)}\abs{\nabla \cbar_{\eta,\tau}}^s)\dx\dt\\  
    &\leq C\int_0^T\big(\norm{\cbar^{(m+\alpha)}_{\eta,\tau}}^s_{L^{\frac{ps}{p-s}}(\Omega)}\norm{\rmD^2\chibar_{\eta,\tau}}^s_{L^p(\Omega;\R^{d\times d\times d})} \\
    &\qquad\qquad + \big\|{\cbar^{\frac{m}{2}}_{\eta,\tau}\big\|}^s_{L^\frac{2s}{2-s}(\Omega)}\big\|{\nabla \cbar_{\eta,\tau}^{\frac{m}{2}}}\big\|^s_{L^2(\Omega)}\big)\dt,
\end{align*}
where in the first and second inequality we used Assumptions \ref{Assu:MobilityTensor} and \ref{Assu:FreeEn},
respectively, and in the third inequality, we used H\"older's inequality. We now consider the two terms on the right-hand
side separately. 
For the first term, the boundedness follows from the condition $0\leq m+\alpha \leq \frac{p-s}{ps}$, 
i.e., $(m+\alpha)\frac{ps}{p-s}\leq 1$ (see Assumption \ref{Assu:FreeEn}(iii)), and the $L^\infty(0,T;L^1(\Omega))$-bound for $\cbar_{\eta,\tau}$ and $L^\infty(0,T;L^p(\Omega;\R^{d\times d\times d}))$-bound for $\rmD^2\chibar_{\eta,\tau}$.  

For the second term, we use the Gagliardo-Nirenberg-Sobolev inequality to bound
\begin{align*}
\big\|{\cbar^{\frac{m}{2}}_{\eta,\tau}\big\|}_{L^\frac{2s}{2-s}(\Omega)} 
&\leq C\big(\big\|{\nabla \cbar_{\eta,\tau}^{\frac{m}{2}}}\big\|^\lambda_{L^2(\Omega)}\big\|{\cbar_{\eta,\tau}^\frac{m}{2}}\big\|^{1-\lambda}_{L^\frac{2}{m}(\Omega)} + \big\|{ \cbar_{\eta,\tau}^{\frac{m}{2}}}\big\|_{L^1(\Omega)}\big)\\
&\leq C\big(1+\big\|{\nabla \cbar_{\eta,\tau}^{\frac{m}{2}}\big\|}^\lambda_{L^2(\Omega)}\big)
\end{align*}
for almost every $t\in [0,T]$. Here, $\lambda\in(0,1)$ is defined via $\frac{2-s}{2s} = \lambda\frac{d-2}{2d} + (1{-}\lambda)\frac{m}{2}$, or using that $s=\frac{md+2}{md+1}$, via $s(1+\lambda)=2$, and we have used that $\cbar_{\eta,\tau}$ is bounded in $L^\infty(0,T;L^1(\Omega))$. 
Thus, using that $1<s<2$, we obtain
\begin{align*}
    \int_0^T \big\|{\cbar^{\frac{m}{2}}_{\eta,\tau}}\big\|^s_{L^\frac{2}{2-s}(\Omega)}\norm{\nabla \cbar_{\eta,\tau}^{\frac{m}{2}}}^s_{L^2(\Omega)}\dt 
    &\leq C\int_0^T \big\|{\nabla \cbar_{\eta,\tau}^{\frac{m}{2}}}\big\|^s_{L^2(\Omega)} +  \big\|{\nabla \cbar_{\eta,\tau}^\frac{m}{2}}\big\|^{s(1+\lambda)}_{L^2(\Omega)} \dt\\
    &\leq C\big(1+\big\|{\nabla \cbar_{\eta,\tau}^\frac{m}{2}}\big\|^{2}_{L^2(0,T;L^2(\Omega))}\big).
\end{align*}
Using the bounds in \eqref{Eqn:TimeDisc:APrioriC}, the bound in \eqref{Eqn:TimeDisc:APrioriMob} follows.

    \underline{\textbf{Case II}: $\gamma_2\geq\gamma_1>0$.} To prove the bound for the flux in 
   \eqref{Eqn:APriori:Reg:mob:CaseII}, we estimate for the exponent $1<s=\min\{\frac{md+2(r+2)}{md+r+2},\frac{d(m+r+1)+2(r+2)}{d(m+r+1)+r+2}\}<2$
\begin{align*}
    &\norm{\calM(\nabla\chibar_{\eta,\tau},\cbar_{\eta,\tau})\nabla\mubar_{\eta,\tau}}_{L^s([0,T]\times\Omega)}^s 
	\leq C\int_0^T\int_\Omega \big(\cbar^m_{\eta,\tau}\big|\pl_{Fc}^2\Phi \rmD^2\chibar_{\eta,\tau} 
	+ \pl_{cc}^2\Phi\nabla \cbar_{\eta,\tau}\big|\big)^s\dx\dt\\
    &\quad\leq C\int_0^T\int_\Omega \big(\cbar_{\eta,\tau}^{s(m+\alpha)}\abs{\rmD^2\chibar_{\eta,\tau}}^s + (\cbar_{\eta,\tau}^{s(m-1)}+\cbar_{\eta,\tau}^{s(m+r)})\abs{\nabla \cbar_{\eta,\tau}}^s)\dx\dt\\  
    &\quad\leq C\int_0^T\big(\norm{\cbar^{(m+\alpha)}_{\eta,\tau}}^s_{L^{\frac{ps}{p-s}}(\Omega)}\norm{\rmD^2\chibar_{\eta,\tau}}^s_{L^p(\Omega;\R^{d\times d\times d})} \\
    &\qquad\qquad + \big\|{\cbar^{\frac{m}{2}}_{\eta,\tau}}\big\|^s_{L^\frac{2s}{2-s}(\Omega)}\big\|{\nabla \cbar_{\eta,\tau}^{\frac{m}{2}}}\big\|^s_{L^2(\Omega)} + \big\|{\cbar^{\frac{m+r+1}{2}}_{\eta,\tau}}\big\|^s_{L^\frac{2s}{2-s}(\Omega)}\big\|{\nabla \cbar^{\frac{m+r+1}{2}}}\big\|^s_{L^2(\Omega)}\big)\dt,
\end{align*}
where in the first and second inequality we used Assumptions \ref{Assu:MobilityTensor} and \ref{Assu:FreeEn}, and in the second and third inequality we used H\"older's inequality. Again, we look at the integrals separately. 
For the first integral, the boundedness follows from the condition $(m+\alpha)\frac{ps}{p-s}\leq 2+r$ (cf. Assumption \ref{Assu:FreeEn}(iii).
For the second integral, we argue as before, now using the bound $m\leq 4+2r$ (which directly follows from $m\leq 3+r$ and $r\geq -1$). Taking $s=\frac{md+2(r+2)}{md+r+2}>1$ the boundedness follows. The boundedness of the third integral follows from $m\leq 3+r$ and $s=\frac{d(m+r+1)+2(r+2)}{d(m+r+1)+r+2}>1$.

\emph{Step 6.} 
 To prove the final bound of \eqref{Eqn:TimeDisc:APrioriC} for the time derivative $\dot \chat_{\eta,\tau}$, we now test \eqref{Eqn:TimeDisc:FiniteStrainChemPot:Reg} with $\psi\in L^{s'}(0,T;H^\theta(\Omega))$ and use Assumption \ref{Assu:Permeability} to find
    \begin{align*}
	\int_0^T\ip{\dot \chat_{\eta,\tau},\psi}\dt 
    &\leq C\norm{\calM(\nabla\chibar_{\eta,\tau},\cbar_{\eta,\tau})\nabla\mubar_{\eta,\tau}}_{L^s([0,T]\times\Omega)}\norm{\psi}_{L^{s'}(0,T;W^{1,s'}(\Omega))} \\
    &\qquad  + 
    C\eta\norm{\mubar_{\eta,\tau}}_{L^2(0,T;H^\theta(\Omega))}\norm{\psi}_{L^2(0,T;H^\theta(\Omega))}\\
    &\qquad + C\big(\norm{\mubar_{\eta,\tau}{-}\mubar_{\tau,\ext}}_{L^2([0,T]\times\pl\Omega)}\Big) \norm{\psi}_{L^2([0,T]\times\pl\Omega)}.
 	\end{align*}	
	Thus, using the previously obtained uniform bounds, we conclude that
	\begin{equation}
	\label{Eqn:FiniteStrain:APrioriCdot}
	\int_0^T \ip{\dot \chat_{\eta,\tau},\psi}\dt \leq C(\norm{\psi}_{L^{s'}(0,T;W^{1,s'}(\Omega))} + \sqrt{\eta}\norm{\psi}_{L^2(0,T;H^\theta(\Omega))}),
	\end{equation}
   for all $\psi\in L^{s'}(0,T;H^{\theta}(\Omega))$ so that $\norm{\dot \cbar_{\eta,\tau}}_{L^{s}(0,T;H^\theta(\Omega)^*)}\leq C$, finishing the proof of \eqref{Eqn:TimeDisc:APrioriC}.
\end{proof}

\section{Limit passage \texorpdfstring{$(\tau,\eta)\to 0$}{Limit passage}}
\label{Sect:LimitPassage}

We now pass to the limit $(\tau,\eta)\to 0$. Note that to obtain strong convergence of the concentration $c_\eta$, we cannot directly use the Aubin--Lions lemma since we do not have a bound for $\nabla \cbar_{\eta,\tau}$, but instead for $\nabla \cbar_{\tau,\eta}^{\frac{m}{2}}$, and, if $\gamma_1>0$, also for $\nabla \cbar_{\eta,\tau}^{\frac{m}{2} + i}$ ($i\in \{\frac{1+r}{2}, 1{+}r\}$). Thus, we use a generalization, Dubinski\u\i's theorem, cf.\ Theorem~\ref{Th:DubinskiiCompactness}.

	\begin{Prop}[Limit passage $(\tau,\eta)\to 0$]
         \label{Prop:LimitTimeDisc}
		There exist subsequences (not relabeled) such that the interpolants  $(\chihat_\tau, \chat_\tau, \mubar_\tau)$ converge to a weak solution $(\chi, c, \mu)$ of \eqref{Eqn:FiniteStrainSystem}--\eqref{Eqn:FiniteStrainInitial} in the sense of Definition \ref{Def:weakSolutionFiniteStrain}).
	\end{Prop}

\begin{proof}
We show the limit passage for Case I, i.e., $\gamma_1=\gamma_2=0$. The limit passage for Case II, i.e., $\gamma_2\geq \gamma_1>0$ follows in a similar way.

\emph{Step 1.}
Using the a priori estimates from Lemma \ref{lem:AprioriEstimatesTimeDisc}, 
we can extract converging subsequences (not relabeled) and some $(\chi,c,\mu,\Theta)$ such that
\begin{alignat*}{2}
			&\chihat_{\eta,\tau} \wstarto \chi\ &&\text{in}\ L^\infty(0,T;W^{2,p}_\id(\Omega;\R^d))\cap H^1(0,T;H^1(\Omega;\R^d)),\\
			&\chat_{\eta,\tau} \wstarto c&&\text{in}\ L^\infty(0,T;L\log L(\Omega))\cap W^{1,s}(0,T;H^\theta(\Omega)^*),\\
			&\mubar_{\eta,\tau} \wto \mu_*&&\text{in}\ L^2([0,T]\times \pl\Omega),\\
            &\calM(\nabla\chibar_{\eta,\tau}, \cbar_{\eta,\tau})\nabla\mubar_{\eta,\tau} \wto \Theta\qquad&&\text{in}\ L^s([0,T]\times \Omega).
		\end{alignat*}
		Using the Aubin--Lions lemma, we can extract a strongly converging subsequence (not relabeled) such that
\[
\chihat_{\eta,\tau} \sto \chi\ \text{in}\ C(0,T;C^{1,\lambda}(\overline\Omega;\R^d)),\quad\text{with }\lambda=1-\frac{d}{p},
\]
Next, we note that since $\cbar_{\eta,\tau}^{\frac{m}{2}}\in L^2(0,T;L^1(\Omega))$ (as $1\leq m\leq 2$), it follows from Poincar\'e's inequality that $\cbar_{\eta,\tau}^{\frac{m}{2}}\in L^2(0,T;H^1(\Omega))$. Thus, we can apply Corollary \ref{Cor:DubinskiiCompactness} of Dubinski\u\i's theorem %
to obtain a strongly converging subsequence (not relabeled)
\[
\chat_{\eta,\tau} \sto c \quad \text{in}\ L^{m}(0,T;L^m(\Omega)).
\]
Note that by Remark 
\ref{Rem:TimeDisc:InterpolantEstimates} these convergences also hold for the piecewise 
constant interpolants $\chibar_{\eta,\tau}$, $\ubar\chi_{\eta,\tau}$, $\cbar_{\eta,\tau}$ and $\ubar c_{\eta,\tau}$.

Consequently, we also have that 
\begin{alignat*}{2}
            &\nabla \cbar_{\eta,\tau}^{\frac{m}{2}} \wto \nabla c^{\frac{m}{2}}\quad &&\text{in}\ L^2(0,T;L^2(\Omega)).
   \end{alignat*}

\medskip

\emph{Step 2.}
We now pass to the limit $\tau\to 0$ in the time-discrete mechanical equation
in \eqref{Eqn:TimeDisc:WeakFiniteStrainStress}. Using the above convergences and the continuity of $\sigma_\el:\mathsf{F}_R\times \R^+\to \R^{d\times d}$
(with $R>0$ as in Step 3 of the proof of the a priori estimates in Lemma~\ref{lem:AprioriEstimatesTimeDisc}), 
the limit in the first term on the left-hand side 
follows from the dominated convergence theorem. Indeed, using the 
$L^\infty(0,T;L^\infty(\Omega;\R^{d\times d}))$-bound for $\nabla\chibar_{\eta,\tau}$ and Assumption \ref{Assu:FreeEn}(iii), we obtain $\abs{\sigma_\el(\nabla\chibar_{\eta,\tau},\cbar_{\eta,\tau})}\leq C(1+\cbar_{\eta,\tau}^{1+\alpha})$,  which is integrable since $0\leq 1+\alpha\leq 1$ (Case I), or since $0\leq 1+\alpha\leq 2+r$ (Case II). 
Using the linearity of $\sigma_\vi(F,\dot F,c)$ 
with respect to $\dot F$, we can also directly pass to the 
limit in the second term. Recall that by Assumption~\ref{Assu:ViscousStress}
the tensor $\widetilde{\mathbb{D}}$ is uniformly bounded, hence, 
the dominated convergence theorem can be applied again. 
The limit in the loading term on the right-hand side follows
directly from the fact that $\wb\ell_\tau\to \ell$ in $L^2(0,T;H^1(\Omega)^*)$. 

Finally, the limit passage in the hyperstress term follows from 
the strong-weak closedness of the convex subdifferential
of $\mathcal{A}:X\to [0,\infty]$ defined by 
\begin{equation*}
\mathcal{A}(\chi) = \begin{cases}
\int_0^T\int_\Omega \mathscr{H}(\rmD^2\chi)\dx\dt\quad &\text{if}\ \chi\in L^\infty(0,T;W^{2,p}(\Omega;\R^d)),\\
+\infty &\text{otherwise},
\end{cases}
\end{equation*}
see Lemma \ref{Lemma:Strong-weakClosednessSubdiff}.
\medskip

\emph{Step 3.}
To pass to the limit $\tau\to 0$ in the diffusion equation
\eqref{Eqn:TimeDisc:WeakFiniteStrainChemPot}, we note that the 
first term on the left-hand side converges 
by the weak convergence of $\dot\chat_{\eta,\tau}$ in $L^2(0,T;H^\theta(\Omega)^*)$. The third term vanishes in the limit $\eta\to 0$ using the first estimate of \eqref{Eqn:TimeDisc:APrioriMu}. The fourth term converges using the weak convergences of $\mubar_{\eta,\tau}$ in $L^2([0,T]\times\pl\Omega)$. Finally, it remains to pass to the limit in the second integral, i.e., to show that
\[
\int_0^T\int_\Omega \calM(\nabla\chibar_{\eta,\tau},\cbar_{\eta,\tau})\nabla\mubar_{\eta,\tau} \cdot \nabla\psi\dx \dt \to \int_0^T\int_\Omega \calM(\nabla\chi,
c)\nabla\mu \cdot \nabla\psi\dx\dt
\]
for all $\psi\in L^{s'}(0,T;H^\theta(\Omega))$. Note that we have (up-to a subsequence) $\calM(\nabla\chibar_{\eta,\tau}, \cbar_{\eta,\tau})\nabla\mubar_\tau \wto \Theta$ in $L^s([0,T]\times\Omega))$ by Step 1.

Using the strong convergences of $\nabla\chibar_{\eta,\tau}$ and $\cbar_{\eta,\tau}$ we can assume that
$\nabla \chibar_{\eta,\tau}(t,x) \to \nabla\chi(t,x)$ in $\GL^+(d)$ and $\cbar_{\eta,\tau}(t,x) \to c(t,x)$ for almost every $(t,x)\in[0,T]\times\Omega$.

To pass to the limit, we now introduce the rescaled mobility tensor $\wt\calM$ by defining $\wt\calM(F,c) = c^m\calM(F,c)$ for all $F\in \mathsf{F}_R$ (some $R>0$) and $c>0$, and the rescaled derivatives of free energy $A$ and $B$ defined by $c^\alpha A(F,c) = \pl^2_{Fc}\Phi(F,c)$ and $B(F,c) = c\pl^2_{cc}\Phi(F,c)$ for all $F\in \GL^+(d)$ and $c>0$. 
In particular, note that by Assumption \ref{Assu:MobilityTensor} we have $\abs{\wt\calM}\leq C$ and by Assumption \ref{Assu:FreeEn} that $\abs{A}, \abs{B} \leq C$ for all $F$ and $c>0$.
Then, it follows that
\begin{multline*}
\int_0^T\int_\Omega \calM(\nabla\chibar_{\eta,\tau},\cbar_{\eta,\tau})\nabla\mubar_{\eta,\tau} \cdot \nabla\psi\dx \dt \\
= \int_0^T\int_\Omega \wt\calM(\nabla\chibar_{\eta,\tau},\cbar_{\eta,\tau})\big(A(\nabla\chibar_{\eta,\tau},\cbar_{\eta,\tau})c_{\eta,\tau}^{m+\alpha}\rmD^2\chibar_{\eta,\tau} + \cbar_{\eta,\tau}^{\frac{m}{2}+1}B(\nabla\chibar_{\eta,\tau},\cbar_{\eta,\tau})\nabla \cbar^{\frac{m}{2}}_{\eta,\tau}\big)\cdot\nabla\psi \dx\dt.
\end{multline*}
Since continuous functions preserve almost everywhere convergence, the limit passage follows by applying Lemma \ref{Lemma:Limitpassage} (therein, stated for the scalar case with the extension to the vectorial case being straightforward) with 
\[
\mathfrak{a}_{\eta,\tau} = \wt\calM(\nabla\chibar_{\eta,\tau},\cbar_{\eta,\tau})A(\nabla\chibar_{\eta,\tau},\cbar_{\eta,\tau})c_{\eta,\tau}^{m+\alpha}\quad \text{and } \mathcal{V}_{\eta,\tau} = \rmD^2\chibar_{\eta,\tau}
\]
and 
\[
\mathfrak{a}_{\eta,\tau} = \wt\calM(\nabla\chibar_{\eta,\tau},\cbar_{\eta,\tau})\cbar_{\eta,\tau}^{\frac{m}{2}+1}B(\nabla\chibar_{\eta,\tau},\cbar_{\eta,\tau})\nabla \cbar^{\frac{m}{2}}_{\eta,\tau}\quad \text{and } \mathcal{V}_{\eta,\tau} = \nabla\cbar_{\eta,\tau}^{\frac{m}{2}}
\]
for the first and second integral, respectively. Note, in particular, that this uses the $L^s([0,T]\times \Omega)$-bound for the flux $\calM(\nabla\chibar_{\eta,\tau},\cbar_{\eta,\tau})\nabla\mubar_{\eta,\tau}$.

The limit passage in the terms arising when $\gamma_2\geq \gamma_1>0$ follows in a similar way, and is therefore omitted.
\medskip

\emph{Step 4.} It remains to show that $\mu_*\in \pl_c\Phi(\nabla\chi,c)$. Note that this is not immediately clear since we do not have an $L^2(0,T;H^1(\Omega))$ estimate for $\mu_\eta$ (as is the case for nondegenerate mobilities), and thus cannot use the trace operator $H^1(\Omega)\hookrightarrow L^2(\pl\Omega)$ to conclude that this $\mu_*$ satisfies $\mu_*\in \pl_c\Phi(\nabla\chi,c)$.

Instead we use that $\nabla\chihat_{\eta,\tau} \sto \chi$ in $C(0,T;C(\wb\Omega;\R^{d\times d}))$ to extract a subsequence (not relabeled) such that $\chihat_{\eta,\tau}(t,x)\to \chi(t,x)$ in $\GL^+(d)$ for almost every $(t,x)\in [0,T]\times \pl\Omega$. Using that $H^1(\Omega)\hookrightarrow L^2(\pl\Omega)$ (see e.g.\ \cite[Theorem 1.1]{Bieg2009ToSFoBoED}), we thus obtain using Dubinski\u\i's theorem a strongly converging subsequence (not relabeled) such that $\cbar_{\eta,\tau}\sto c$ in $L^m([0,T]\times \pl\Omega)$. In particular, we can extract another subsequence sucht that $\cbar_{\eta,\tau}(t,x)\to c(t,x)$ for almost every $(t,x)\in [0,T]\times \pl\Omega$. 
Using the Banach--Saks theorem, we now extract a subsequence $(\eta_n,\tau_n)$ such that $\frac{1}{N}\sum_{n=1}^N \mu_{{\eta_n},{\tau_n}}\sto \mu_*$ in $L^2([0,T]\times \pl\Omega)$. Again extracting a converging (non-relabeled) subsequence, 
we thus find that $\frac{1}{N}\sum_{n=1}^N \mu_{{\eta_n},{\tau_n}}(t,x)\to \mu_*(t,x)$ for all $(t,x)\in [0,T]\times \pl\Omega$. Since continuity of $\pl_c\Phi$ and the almost everywhere convergence of $\nabla\chibar_{\eta,\tau}$ and $\cbar_{\eta,\tau}$ also imply that $\mu_{\eta,\tau} = \pl_c\Phi(\nabla\chibar_{\eta,\tau},\cbar_{\eta,\tau}) \to \pl_c\Phi(\nabla\chi,c)$ for almost every $(t,x)\in [0,T]\times \pl\Omega$, we thus conclude that $\mu_*\in \pl_c\Phi(\nabla\chi,c)$.

\end{proof}

\appendix
\section{Tools}

\begin{Lemma}[Strong-weak closedness of the convex subdifferential]
	\label{Lemma:Strong-weakClosednessSubdiff}
	Let $X$ be a reflexive Banach space, and let $\mathcal{A}:X\to \R\cup\infty$ be convex and lower semicontinuous. If $u_n\sto u$ in $X$ and $\xi_n \in \pl \mathcal{A}(u_n)\subset X^*$ such that $\xi_n\wto \xi$ in $X^*$, then $\xi\in \pl \mathcal{A}(u)$.
	\end{Lemma}	
	\begin{proof}
	Using convexity of $\mathcal{A}$, we obtain $\mathcal{A}(v)\geq \mathcal{A}(u_n) + \ip{\xi_n, v-u_n}$ for all $v\in X^*$. The result now follows using by taking the $\liminf$ of this inequality, and using the lower semicontinuity and given convergences.
	\end{proof}

Here we state Dubinski\u\i's theorem \cite{Dubi1965WCfNEaPE} in the form of \cite{BarSul2012RoDNCET}.
\begin{Def}
Let $A$ be a Banach space. Then, $\calM_+$ is a \textbf{seminormed nonnegative cone} in $A$ if $\calM_+\subseteq A$ satisfies:
\begin{itemize}
\item[(i)] For all $u\in\calM_+$, $\gamma\geq 0$, we have $\gamma u\in \calM_+$.
\item[(ii)] There exists a function $[\cdot]:\calM_+\to [0,\infty)$ (the semi-norm) such that $[u]=0$ if and only if $u=0$, and for $\gamma\geq 0$ we have $[\gamma u] = \gamma[u]$.
\end{itemize}
We say that $\calM_+\hookrightarrow A$ continuously if there exists $C>0$ such that $\norm{u}_A\leq C[u]$ for all $u\in \calM_+$, and we say that $\calM_+\hookrightarrow A$ compactly if every bounded sequence in $\calM_+$ has a converging subsequence in $A$.
\end{Def}
\begin{Th}[Dubinski\u\i] Let $\mathcal{M}_+$ be a seminormed nonnegative cone, and $A_0$, $A_1$ Banach spaces such that $\mathcal{M}_+ \hookrightarrow A_0$ compactly, and $A_0\hookrightarrow A_1$ continuously. Let 
\[
\mathcal{Y}_+ = \{ \varphi:[0,T]\to \mathcal{M}_+ \mid [\varphi]_{\mathcal{Y}_+} := [\varphi]_{L^p(0,T;\mathcal{M}_+)} + \norm{\dot\varphi}_{L^{p_1}(0,T;A_1)} \leq C\}
\]
for $1\leq p$, $p_1\leq \infty$. Then, $\mathcal{Y}_+$ is a seminormed nonnegative cone in $L^p(0,T;A_0)\cap W^{1,p_1}(0,T;A_1)$, and $\mathcal{Y}_+ \hookrightarrow L^p(0,T;A_0)$ compactly if either $p_1>1$ or $p<\infty$.
\label{Th:DubinskiiCompactness}
\end{Th}

\begin{Cor}
Let $\omega\geq \frac{1}{2}$, $\theta>\frac{d}{2}$, and let $(\varphi_k)$ be a sequence of nonnegative functions on $[0,T]\times\Omega$ such that there exists a constant $C>0$ for which 
\[
\norm{(\varphi_k)^\omega}_{L^2(0,T;H^1(\Omega))} + \norm{\dot\varphi_k}_{L^2(0,T;H^\theta(\Omega)^*)} \leq C.
\]
Then, $(\varphi_k)$ is relatively compact in $L^{2\omega}(0,T;L^{2\omega}(\Omega))$.
\label{Cor:DubinskiiCompactness}
\end{Cor}   
\begin{proof}
The proof is given in the time-discrete case in \cite{ChJuLi2014NoALDl}. In our case, we just apply Theorem \ref{Th:DubinskiiCompactness} with $\mathcal{M}_+ = \{\varphi \mid \varphi^\omega \in H^1(\Omega)\}$, $A_0 = L^{2\omega}(\Omega)$, $A_1 = H^\theta(\Omega)^*$, and the seminorm 
\[
[\varphi]_{\mathcal{Y}_+} = \norm{\varphi^\omega}^{\frac{1}{\omega}}_{L^2(0,T;H^1(\Omega))} + \norm{\dot\varphi}_{L^2(0,T;H^\theta(\Omega)^*)}.
\]
\end{proof}

\begin{Lemma}[{\cite[Lemma~A.3]{FHKM2022GEAE}}] Assume that $\mathfrak{a}_n, \mathcal{V}_n: Q\to \R$, $n\in\N$ 
are measurable functions in $Q\subset\R^m$ (open, bounded)
such that $\mathcal{V}_n \wto \mathcal{V}$ in $L^1(Q)$, $\mathfrak{a}_n(x) \to \mathfrak{a}(x)$
for almost every $x\in Q$, and $\sup_{n\in\N}\|\mathfrak{a}_n\mathcal{V}_n\|_{L^{1+\sigma}(Q)}<\infty$ for some $\sigma>0$. 
Then $\mathfrak{a}_n\mathcal{V}_n 
\wto \mathfrak{a} \mathcal{V}$ in $L^{1+\sigma}(Q)$.
\label{Lemma:Limitpassage}
\end{Lemma}

\subsection*{Acknowledgement}
M.L. was partially supported by DFG via the Priority Program SPP 2256 \emph{Variational Methods for Predicting Complex
Phenomena in Engineering Structures and Materials} (project no. 441470105, subproject Mi 459/9-1 \emph{Analysis for thermomechanical models with internal variables}).

\printbibliography[]
\end{document}